\documentclass[final,review]{elsarticle}
\usepackage{amsmath}
\usepackage{amssymb}
\usepackage{latexsym}
\usepackage{graphicx}
\usepackage{psfrag}
\usepackage{subfigure}
\usepackage{color}
\usepackage{bm}

\usepackage{showkeys}

\usepackage{paralist}
\usepackage[colorlinks=true]{hyperref}

 \newtheorem{thm}{Theorem}[section]

 \newtheorem{rem}[thm]{Remark}
 
\def\bbf{\mathbf{f}}
\def\bg{\mathbf{g}}

\def\bu{\mathbf{u}}

\def\bL{\mathbf{L}}

\def\bR{\mathbf{R}}

\def\n{\noindent}
\def\pt{\partial}

\def\Re{\mathbb{R}}

\def\f#1#2{\frac {#1}{#2}}
\def\f32{\frac 32}

\def\d{\displaystyle}

\def\bga{\begin{array}}
\def\eda{\end{array}}

\def\Gm{\Gamma}
\def\gm{\gamma}

\def\la{\lambda}

\def\al{\alpha}
\def\td{\tilde}

\def\d{\displaystyle}
\def\dfr#1#2{\displaystyle{\frac{#1}{#2}}}

\begin{document}
\begin{frontmatter}
\title{A two-stage fourth order time-accurate discretization for Lax--Wendroff type flow solvers\\ II. High order numerical boundary conditions}

\author[bnu]{Zhifang Du}
\ead{du@mail.bnu.edu.cn}
\author[iapcm,capt]{Jiequan Li\corref{cor1}}
\ead{li_jiequan@iapcm.ac.cn}

\cortext[cor1]{Corresponding author}

\address[bnu]{School of Mathematical Sciences, Beijing Normal University, 100875, Beijing, P. R. China}
\address[iapcm]{Laboratory of Computational Physics,  Institute of Applied Physics and Computational Mathematics, Beijing, P. R. China}
\address[capt]{Center for Applied Physics and Technology, Peking University, Beijing, P. R. China}

\begin{abstract}
This paper serves  to treat  boundary conditions numerically with high order accuracy  in order to suit the  two-stage fourth-order finite volume schemes for hyperbolic problems  developed in [{\em J. Li and Z. Du, A two-stage fourth order time-accurate discretization {L}ax--{W}endroff type flow solvers, {I}. {H}yperbolic conservation laws, SIAM, J. Sci. Comput., 38 (2016), pp.~A3046--A3069}]. As such,  it is significant when capturing small scale structures near physical boundaries. Different from previous contributions in literature,  the current approach constructs a fourth order accurate approximation to boundary conditions  by only using the Jacobian matrix of the flux function (characteristic information) instead of its successive differentiation of governing equations leading to tensors of high ranks in the inverse Lax-Wendroff method. Technically, data in several ghost cells are constructed with interpolation so that the interior scheme  can be implemented over  boundary cells, and theoretical boundary condition has to be modified properly at intermediate stages so as to make the  two-stage scheme over  boundary cells  fully consistent with that over interior cells.  This is nonintuitive and   highlights the fact that \em theoretical  boundary conditions are only prescribed for continuous partial differential equations (PDEs), while they must be approximated in a consistent way (even though they could be exactly valued) when the PDEs are discretized.  Several numerical examples are provided to illustrate the performance of the current approach when dealing with general boundary conditions.
\end{abstract}

\end{frontmatter}

\section{Introduction}\label{sec:intro}

This study aims to develop high order accurate numerical approximations to theoretical boundary conditions suitable for the two-stage fourth-order accurate finite volume scheme based on the Lax--Wendroff type flow solvers in \cite{Du-Li-1}. This scheme differs from Runge--Kutta type schemes  in the sense that   second-order Lax--Wendroff type flow solvers \cite{Du-Li-1, Xu-Li}  are used to  achieve  fourth-order  temporal accuracy through a two-stage temporal iteration. The numerical evidence shows that it is less computational consuming, less dissipative near discontinuities and more effective to resolve physical and multi-dimensional effects.
Yet, the numerical boundary conditions have not been constructed correspondingly and  in particular,  the theoretical boundary data have to be modified at intermediate stages,  analogous to other multi-stage methods \cite{BD-ACC, BD-ACC-2}, through our detail analysis.
\\

\begin{figure}[!htb]
\centering
 \includegraphics[scale=0.5]{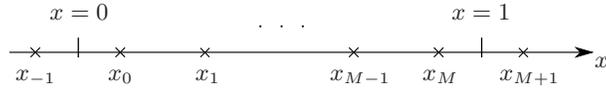}
\small
    \put(-206,15){$ x=0 $}
    \put(-219,-8){$ x_{-1} $}
    \put(-184,-8){$ x_{0} $}
    \put(-151,-8){$ x_{1} $}
    \put(-54,15){$ x=1 $}
    \put(-100,-8){$ x_{M-1} $}
    \put(-65,-8){$ x_{M} $}
    \put(-36,-8){$ x_{M+1} $}
    \put(0,-3){$x$}
 \caption{The computational domain $(0,1)$. Set $ x_0=h/2 $ and $ x_M=1-h/2 $. Then $ x_{-1}=-h/2 $ and $ x_{M+1}=1+h/2 $ are ghost points.}
 \label{fig:kreiss}
\end{figure}

There are a large number  of contributions to numerical boundary conditions for hyperbolic problems  in literature and most of them are only first- or second-order accurate. For example, second-order accurate finite difference approximations are constructed  for  the wave equation with Dirichlet or Neumann boundary conditions in \cite{EB-wave-Dirichlet, EB-wave-Neumann} from which our basic idea comes. As is well-known, a second order wave equation is equivalent to a system of first order hyperbolic equations. Therefore  the similar idea applies to approximate (theoretical) boundary conditions  of  hyperbolic conservation laws in \cite{EB-consv-FD}. We briefly illustrate their idea in the finite difference framework by considering  the initial boundary value problem (IBVP) for a scalar conservation law
\begin{equation}
\label{eq:IBVP}
\left\{
\begin{array}{ll}
  \dfr{\pt u}{\pt t} + \dfr{\pt f(u)}{\pt x}=0,\ \  \ &  x \in (0,1),~ t > 0,\\[3mm]
  u(x,0) = u_0(x), & x \in (0,1),\\[3mm]
  u(0,t) = g(t), & t > 0.
\end{array}
\right.
\end{equation}
Assume that $ f^\prime(u)>0 $ for all $ u\in\Re $ so  that $ x=0 $ is an inflow boundary and $ x=1 $ is  an outflow boundary. We  equally distribute  $ M+1 $ points $\{x_j = (j+1/2)h : j=0,1,\dots,M\}$ in the computational domain $(0,1)$, as shown in Figure \ref{fig:kreiss}. We use $u_j$ to denote the value of $u$ at $x=x_j$ and suppress the index for the time levels. Obviously
at the inflow boundary, the solution value at the ghost point $ x_{-1} $ is required in order  to perform a second-order finite difference at $ x_0 $. For this purpose, a polynomial is constructed in the region around the inflow boundary by using point-wise values $ u_{-1} $, $ u_0 $ and $ u_1 $,
\begin{equation}
L(x) = g_{-1}(x)u_{-1} + g_0(x)u_0 + g_1(x)u_1,
\end{equation}
from which we want to find  the value $ u_{-1} $. The Lagrangian interpolation tells that
\begin{equation}
\begin{array}{l}
g_{-1}(x) = \dfr{(x-x_0)(x-x_1)}{2h^2}, \\[3mm]  g_0(x) = \dfr{(x-x_{-1})(x-x_1)}{-h^2}, \\[3mm]  g_1(x) = \dfr{(x-x_{-1})(x-x_0)}{2h^2}.
\end{array}
\end{equation}
Then $ u_{-1} $ can be obtained by solving the linear equation $  u(0,t) = L(0) $ where  $ u(0,t) = g(t)$. At the outflow boundary $ x=1 $, we simply use the extrapolation
\begin{equation}
  u_{M+1} = 2u_M - u_{M-1}
\end{equation}
 to obtain the value $ u_{M+1} $ since the signal goes out of the computational domain at this end.  Other works can be found   e.g. in \cite{EB-hbox, EB-hbox-md, EB-uncut-FV}  on first- and second-order accurate approaches for hyperbolic conservation laws in the finite volume framework  and in \cite{BF-DG} on the body fitted boundary treatment for the discontinuous Galerkin (DG) methods. \\

In the present paper, we propose a fourth-order accurate boundary condition treatment in order to suit the two-stage fourth-order finite volume scheme developed in \cite{Du-Li-1}. This boundary condition treatment can be regarded as a fourth-order extension of the methods developed in \cite{EB-wave-Dirichlet, EB-wave-Neumann, EB-consv-FD}. With the methodology described above, a polynomial is constructed near the boundary to interpolate the data in ghost cells with the information from both the interior cells and the boundary condition $ u(0,t)=g(t) $. In order to reduce the number of the interior cells used,  first-order spatial derivatives of the conservative variables are utilized in addition to  their values at the boundary, for which we adopt the inverse Lax-Wendroff (ILW) approach in \cite{ILW-HJ-1,ILW,high-order-ILW,ILW-stb}. Recall that the original ILW method makes successive differentiation of the governing equations in order to calculate high order derivatives of the solution.  A novelty of the current approach is that just first order derivatives are used in the approximation and governing equations themselves are  used only once in the process of the substitution of spatial derivatives by temporal derivatives.   There is no need to make the complicated successive differentiation procedure. 

Another novelty is about the treatment of theoretical boundary data when approximated as  numerical boundary conditions.  In analogy with other multi-stage methods \cite{BD-ACC,BD-ACC-2}, the theoretical boundary data  needs  modifying properly at the intermediate stage for the two-stage approach in \cite{Du-Li-1}  so that fourth order accuracy is achieved. Otherwise, the direct input of the boundary data would lose accuracy by our detailed analysis.  This seems quite counter-intuitive. However, this highlights the fact that {\em theoretical boundary conditions are prescribed for continuous partial differential equations (PDEs). As the  
PDEs are discretized for numerical purposes, the boundary data should be properly approximated, instead that   their exact values are directly input.}    In other words, the approximation of theoretical boundary data should be consistent with the discretization of the associated governing equations.

This paper is organized as follows. After the introduction section, we will have a brief review over the two-stage fourth-order scheme for hyperbolic conservation laws in Section \ref{sec:scheme}. In  Section \ref{sec:scalar}, the basic idea and formulation are explained using  the  one-dimensional scalar equation. The system case is shown in Section \ref{sec:system}, with the specification to the compressible Euler equations.  Several numerical examples are shown in Section \ref{sec:numer} to display the performance of the current numerical boundary condition treatment.

\vspace{2mm}

\section{A review over the two-stage fourth-order scheme}\label{sec:scheme}

In \cite{Du-Li-1},  we proposed a two-stage fourth order accurate temporal discretization based on Lax--Wendroff type flow solvers, particularly for hyperbolic conservation laws,
\begin{equation}\label{eq:govern-eq}
\begin{array}{l}
  \dfr{\pt \bu}{\pt t} + \dfr{\pt \bbf(\bu)}{\pt x}=0, \ \ \  x\in \mathbb{R},  t > 0,\\[3mm]
  \bu(x,0) =\bu_0(x), \ \ \ \ x\in \Re,
  \end{array}
\end{equation}
where $\bu$ is  a conservative variable and $\bbf(\bu)$ is the associated  flux function vector. Given a computational mesh $I_j=(x_{j-\frac 12}, x_{j+\frac 12})$ with 
the size $h =x_{j+\frac 12}-x_{j-\frac 12}$, we write \eqref{eq:govern-eq} in the semi-discrete form
\begin{equation}
\dfr{d\bar \bu_j(t)}{dt} =\mathcal{L}_j(\bu):= -\dfr{1}{h}\Big[\bbf(\bu(x_{j+\frac 12},t)) - \bbf(\bu(x_{j-\frac 12},t))\Big],
\label{scheme}
\end{equation}
where  $\bar\bu_j(t)$  is the average of the solution $\bu(x,t)$ over the control volume $I_j$,
\begin{equation}
\bar\bu_j(t)=\dfr{1}{h}\int_{I_j} \bu(x,t)dt,
\end{equation}
and  $\bu(x_{j+\frac 12},t)$ is the solution of \eqref{scheme} in a certain sense.  Thanks  to the Lax--Wendroff flow solvers, which is specialized as the GRP solver \cite{Li-1} in the current study, $\mathcal{L}_j(\bu)$ is well-defined. We denote by $t^n=nk$, $n=1,2,\cdots$, and by $k$ the  time increment. Then the two-stage approach for \eqref{eq:govern-eq} is summarized as follows.

\vspace{2mm}

\n {\bf Step 1.}  With cell averages $\bar\bu^n_j$ and interface values $\hat\bu^n_{j+\frac 12}$ at the time level $t=t^n$, reconstruct the initial data as a piece-wise polynomial $\bu(x,t^n) =\bu^n(x)$ using the HWENO interpolation \cite{du-HWENO}. Then we compute instantaneous Riemann solutions $\bu_{j+\frac 12}^{n}$ and time derivatives $(\pt\bu/\pt t)_{j+\frac 12}^{n}$ using the GRP solver \cite{Li-1}.
\vspace{0.2cm}

\n {\bf Step 2.} Compute intermediate cell averages $\bar{\bu}_j^{n+\frac 12}$ and interface values $\hat\bu_{j+\frac 12}^{n+\frac 12}$ at $t^{n+\frac 12}=t^n+\frac k2 $ using the following formulae,
\begin{equation}
\begin{array}{l}
\bar\bu_j^{n+\frac 12} =\bar\bu_j^n -\dfr{ k }{2h}\Big[\bbf_{j+\frac 12}^{(1)}-\bbf_{j-\frac 12}^{(1)}\Big],\\[3mm]
\d \bbf_{j+\frac 12}^{(1)} = \bbf(\bu_{j+\frac 12}^{n}) + \dfr{ k }{4} \dfr{\pt\bbf}{\pt\bu}({\bu}^{n}_{j+\frac 12})\left(\dfr{\pt\bu}{\pt t}\right)^{n}_{j+\frac 12},\\[3mm]
\hat\bu_{j+\frac 12}^{n+\frac 12} = \bu_j^{n} + \dfr{k}{2}\left(\dfr{\pt\bu}{\pt t}\right)^{n}_{j+\frac 12},
\end{array}
\label{eq:1d-al}
\end{equation}
where $k$ is the time increment constrained by the CFL condition, and ${\pt\bbf}/{\pt\bu}$ is the Jacobian matrix of $\bbf(\bu)$. With cell averages $\bar{\bu}_j^{n+\frac 12}$ and interface values $\hat\bu_{j+\frac 12}^{n+\frac 12}$, we use the HWENO interpolation again to construct a piece-wise polynomial $\bu^{n+\frac 12}(x)$ and continue to find Riemann solutions $\bu_{j+\frac 12}^{n+\frac 12}$ and time derivatives $(\pt\bu/\pt t)_{j+\frac 12}^{n+\frac 12}$ at the intermediate stage $t^{n+\frac 12}=t^n+\frac{ k }2$, as done in Step 1.
\vspace{0.2cm}

\n {\bf Step 3.}  Advance the solution to the next time level $t^{n+1}= t^n+ k $ by
\begin{equation}
\begin{array}{l}
\bar\bu^{n+1}_j =\bar \bu_j^n -\dfr{k}{h}\Big[\bbf_{j+\frac 12}^{4th} -\bbf_{j-\frac 12}^{4th}\Big],\\[3mm]
\d \bbf_{j+\frac 12}^{4th}=\bbf(\bu_{j+\frac 12}^{n}) + \dfr{k}{2}\left[\dfr 13 \dfr{\pt\bbf}{\pt\bu}({\bu}^{n}_{j+\frac 12})(\dfr{\pt\bu}{\pt t})^{n}_{j+\frac 12} + \dfr 23 \dfr{\pt\bbf}{\pt\bu}({\bu}^{n+\frac 12}_{j+\frac 12})(\dfr{\pt\bu}{\pt t})^{n+\frac 12}_{j+\frac 12}\right], \\[3mm]
\hat\bu_{j+\frac 12}^{n+1} = \bu_j^{n} + k\left(\dfr{\pt\bu}{\pt t}\right)^{n+\frac 12}_{j+\frac 12}.
\end{array}
\label{eq:1d-al2}
\end{equation}

\vspace{2mm}

\begin{rem}
The HWENO reconstruction, based on the GRP solver,  was  developed  in \cite{du-HWENO}, to which is referred for details.
In the smooth case, given cell averages $ \bar{\bu}_j $ and $ x $-differences $ \Delta\bu_{j} $ of the solution,
\begin{equation}
\begin{array}{ll}
  \bar{\bu}_{j} = \dfr{1}{h_j}\displaystyle\int_{I_{j}}\bu(x, t)dx, \ \Delta\bu_{j} = \dfr 1h \displaystyle\int_{I_{j}}\bu_x(x, t)dx = \dfr 1h \left[\hat\bu(x_{j+\frac 12}, t) - \hat\bu(x_{j-\frac 12}, t)\right],
\end{array}
\end{equation}
cell boundary values are reconstructed as
\begin{equation}\label{eq:HWENO-general}
\bga{l}
  \bu_{j-\frac 12,-} = \dfr{1}{120}(-23\bar{\bu}_{j-2}+76\bar{\bu}_{j-1}+67\bar{\bu}_{j}-9h\Delta\bu_{j-2}-21h\Delta\bu_{j}),\\[2mm]
  \bu_{j-\frac 12,+} = \dfr{1}{120}(67\bar{\bu}_{j-1}+76\bar{\bu}_{j}-23\bar{\bu}_{j+1}+21h\Delta \bu_{j-1}+9h\Delta\bu_{j+1}),\\[2mm]
  (\pt\bu/\pt x)_{j-\frac 12,\pm} = \dfr{1}{12h}(\bar{\bu}_{j-2} - 15\bar{\bu}_{j-1}+15\bar{\bu}_{j}-\bar{\bu}_{j+1}).
\eda
\end{equation}
 In the presence of discontinuities, a WENO-type stencil selection is performed. Details can be found in \cite{du-HWENO}.

\end{rem}
\vspace{2mm}

\begin{rem}
Riemann solutions $\bu_{j+\frac 12}^{n}$, and time derivatives $(\pt\bu/\pt t)_{j+\frac 12}^{n}$ in Step 1 are defined as
\begin{equation}\label{eq:GRP-solution-0}
\bga{l}
\bu^{n}_{j+\frac 12} = \d\lim_{t\rightarrow t^n+0} \bu(x_{j+\frac 12},t),\\[3mm]
\left(\dfr{\pt\bu}{\pt t}\right)^{n}_{j+\frac 12} = \d\lim_{t\rightarrow t^n+0} \dfr{\pt\bu}{\pt t}(x_{j+\frac 12},t)= \d\lim_{t\rightarrow t^n+0} - \dfr{\pt\bbf}{\pt\bu}(\bu(x_{j+\frac 12},t))\dfr{\pt \bu}{\pt x}(x_{j+\frac 12},t).
\eda
\end{equation}
The last equality explains the Lax-Wendroff methodology  for which  the governing equation \eqref{eq:govern-eq} is employed.
In practice, due to the singularity of the initial data  $\bu^n(x)$ at $ (x_{j+\frac 12},t^n) $, the GRP solver is applied to resolve the singularity so that  both values
$\bu_{j+\frac 12}^{n+\frac 12}$ and $(\pt\bu/\pt t)_{j+\frac 12}^{n+\frac 12}$ are defined.
\end{rem}

In \cite{Du-Li-1},  we did not touch the  indispensable numerical boundary conditions. In the coming sections, we will present numerical boundary conditions suitable for the implementation of this two-stage scheme.

\vspace{2mm}
\section{Numerical boundary condition treatment for one-dimensional scalar conservation laws}\label{sec:scalar}
Let's consider \eqref{eq:govern-eq} in the computational domain $ (0,1) $ and assume that $\bbf(\bu)$ to be scalar $\bbf(\bu)=f(u)$ with   $ f^{\prime}(u)  > 0$ for all $ u\in \Re $.  Then   $ x=0 $ is an inflow boundary and $ x=1 $  is an outflow boundary. Thus a boundary condition is required at $ x=0 $.  Denote  by $ u(0,t)=g(t) $ the boundary condition for $ t>0 $. The resulting initial boundary value problem (IBVP)  can be formulated as in \eqref{eq:IBVP}.
The mesh $ \{I_j=(x_{j-\frac 12}, x_{j+\frac 12}): x_{j-\frac 12} = jh, j = 0, 1, \dots, M\} $ is used in our computation. Notice that in order to place the left boundary at $ x=0 $, our mesh is a little different from $ I_j = ((j-\frac 12)h, (j+\frac 12)h) $ which is conventionally used in the finite volume context.

In computations, we need numerical boundary conditions at the inflow boundary $ x=0 $ corresponding to the boundary condition $ \bu(0,t) = g(t), t > 0 $. In this section, we will make the numerical treatment for inflow boundary conditions and present the modification at intermediate stages. The WENO-type extrapolation is introduced  to deal with outflow boundary conditions in the last part of this section, which is the same as that in \cite{high-order-ILW}.

\vspace{2mm}

\subsection{Inflow boundary condition treatment}\label{sec:inflow}
As far as the IBVP \eqref{eq:IBVP} is concerned with, several values outside of the computational domain are needed. For the present case,  $ \bar{u}_{-1} $, $ \bar{u}_{-2} $, $ \Delta u_{-1} $ and $\Delta u_{-2}$, defined over  $I_{-1}=(-h,0)$ and $I_{-2}=(-2h,-h)$ respectively,  are needed in the reconstruction procedure \eqref{eq:HWENO-general} for the values indexed by $ j = 0 $ and $ j = 1 $. We call $I_{-1}$ and $I_{-2}$ ghost cells. In this subsection, we suppress all superscripts for the time levels. \\

First of all, assume the solution to be  smooth near the boundary $x=0$. To obtain the values mentioned above, a cubic polynomial,
\begin{equation}\label{eq:poly-5}
 p(x) = \al_3 x^3 +\al_2 x^2 +\al_1 x + \al_0,
\end{equation}
is constructed over $ I_{-2} \cup I_{-1} \cup I_0 \cup I_1 = (-2h,2h) $ to interpolate the solution $ u(x,t^n) $
such that
\begin{equation}\label{eq:condition-5}
 \dfr{1}{h}\displaystyle\int_{I_{i}}p(x)dx = \bar{u}_{i}, \ \ \ \ i=-2, -1, 0, 1.
\end{equation}
Substitute the constraints \eqref{eq:condition-5} into \eqref{eq:poly-5} to determine the coefficients  $ \al_0 $, $ \al_1 $, $ \al_2 $ and $ \al_3 $  as
\begin{equation}\label{eq:poly-5-coeff-0}
\begin{array}{ll}
\alpha_3 = \dfr{\bar{u}_1 - 3\bar{u}_0 + 3\bar{u}_{-1} - \bar{u}_{-2}}{6h^3}, \ \ &\alpha_2 = \dfr{\bar{u}_1 - \bar{u}_0 - \bar{u}_{-1} + \bar{u}_{-2}}{4h^2},\\[3mm]
\alpha_1 = \dfr{-\bar{u}_1 + 15\bar{u}_0 - 15\bar{u}_{-1} + \bar{u}_{-2}}{12h}, \ \ &\alpha_0 = \dfr{-\bar{u}_1 + 7\bar{u}_0  + 7\bar{u}_{-1}  -\bar{u}_{-2}}{12},
\end{array}
\end{equation}
in which $ \bar{u}_{-1} $ and $ \bar{u}_{-2} $ are yet to be determined and they are obtained by evaluating $ p(0) $ and $ p^\prime(0) $ at the boundary $x=0$, 
\begin{equation}\label{eq:boundary-interpolation}
\begin{array}{l}
  p(0) = \dfr{1}{12}(-\bar{u}_{1} + 7\bar{u}_{0} + 7\bar{u}_{-1} - \bar{u}_{-2}) = g(t)+\mathcal{O}(h^4),\\[2mm]
  p^\prime(0) = \dfr{1}{12h}(-\bar{u}_{1}+15\bar{u}_{0}-15\bar{u}_{-1}+\bar{u}_{-2})= -f^\prime(g(t))^{-1} \ g^{\prime}(t)+\mathcal{O}(h^3).
\end{array}
\end{equation}
The first equation results directly from the boundary condition $ u(0,t) = g(t) $ and the interpolation accuracy  $ p(0) = u(0,t) + \mathcal{O}(h^4) $;  the second one is obtained by using the inverse Lax-Wendroff approach \cite{ILW} that expresses  the spatial variation through the temporal variation. The second equation in \eqref{eq:boundary-interpolation} is well-defined by recalling that  $f^{\prime}(u)>0$ is assumed for all $u\in\Re$ at the beginning of this section.
Solving \eqref{eq:boundary-interpolation} in terms of $\bar{u}_{-1}$ and $\bar{u}_{-2}$ yields (by ignoring high order terms)
\begin{equation}\label{eq:ILW-result-u}
\begin{array}{l}
  \bar{u}_{-1} = \dfr 14 (-6g + 6~h~ f'(g)^{-1} \ g^{\prime} + 11\bar{u}_{0} - \bar{u}_{1}),\\[3mm]
  \bar{u}_{-2} = \dfr 14 (-90g + 42~h~ f'(g)^{-1} \ g^{\prime} + 105\bar{u}_{0} - 11\bar{u}_{1}).
\end{array}
\end{equation}
Substituting \eqref{eq:ILW-result-u} into \eqref{eq:poly-5-coeff-0}, in turn,  gives us the explicit expressions of $ \alpha_i$, $i=0,\dots, 3$, and then the expression of $ p(x) $. Therefore we have (by ignoring high order terms)
\begin{equation}\label{eq:ILW-result-ux}
\bga{l}
  \Delta u_{-1} = \dfr{p(0) - p(-h)}h\\[2mm]
\ \ \ \ \ \ \ \ = \dfr 1{8h} (66g - 34~h~ f'(g)^{-1} \ g^{\prime} - 73\bar{u}_{0} + 7\bar{u}_{1}),\\[2mm]
  \Delta u_{-2} = \dfr{p(-h) - p(-2h)}h\\[2mm]
\ \ \ \ \ \ \ \ = \dfr 1{8h} (294g - 118~h~ f'(g)^{-1} \ g^{\prime} - 331\bar{u}_{0} + 37\bar{u}_{1}).
\eda
\end{equation}
The errors of the above approximations are of order  $ \mathcal{O}(h^3) $ due to the accuracy of $ p(x) $. Thus \eqref{eq:ILW-result-u} and \eqref{eq:ILW-result-ux} together provide the  values in the ghost cells $I_{-1}$ and $I_{-2}$.\\
\vspace{2mm}

As there are discontinuities close to the inflow boundary, a WENO-type stencil selecting procedure can be applied. Assume that there is a discontinuity in either $ I_0 $ or $ I_1 $, we  shorten the stencil cell by cell. Denote the stencils by
\begin{equation}
\bga{l}
S^{(2)} = \{I_{-2}, I_{-1}, I_0, I_1\},\ \
S^{(1)} = \{I_{-2}, I_{-1}, I_0\},\ \
S^{(0)} = \{I_{-2}, I_{-1}\}.
\eda
\label{stencil}
\end{equation}
Denote by $ p^{(r)} (x)$  the interpolation polynomial  on  $ S^{(r)} $, $r=0,1,2$, just as the polynomial $ p(x) $ constructed before.   Then define
\begin{equation}
\bga{ll}
  \bar{ u}_{-1}^{(r)} = \dfr{1}{h}\displaystyle\int_{I_{-1}}p^{(r)}(x)dx, & \bar{ u}_{-2}^{(r)} = \dfr{1}{h}\displaystyle\int_{I_{-2}}p^{(r)}(x)dx,\\[3mm]
  \Delta u_{-1}^{(r)} = \dfr 1h (p^{(r)}(0) - p^{(r)}(-h)), & \Delta u_{-2}^{(r)} = \dfr 1h (p^{(r)}(-h) - p^{(r)}(-2h)).
\eda
\end{equation}
The expressions of $ \bar{ u}_{-1}^{(r)} $, $ \bar{ u}_{-2}^{(r)} $, $ \Delta  u^{(r)}_{-1} $ and $ \Delta  u^{(r)}_{-2} $ for $ r=0,1,2 $ will be listed in \ref{app:low-order}.

The smoothness indicators are defined in the same way as for the classical WENO interpolation,
\begin{equation}
\label{eq:SI}
\beta^{(r)} = \displaystyle \sum_{l=1}^{r+1} \int_{I_{-1}} h^{2l-1} \left(\dfr{d^l p^{(r)}}{dx^l}\right)^2 dx, \ \ r=0,1,2,
\end{equation}
where $p^{(r)}$ is the interpolation polynomial of degree $r+1$ on stencil $S^{(r)}$. Expressions of $ \beta^{(r)} $ are put  in \ref{app:SI}. With the Taylor expansion, we can express them at the boundary as
\begin{equation}
\bga{l}
\beta^{(2)} = h^2 {(\dfr{\pt u}{\pt x})}^2 - h^3\dfr{\pt u}{\pt x}\dfr{\pt^2 u}{\pt x^2} + h^4\Big[\dfr 43 {(\dfr{\pt^2 u}{\pt x^2})}^2 + \dfr 13\dfr{\pt u}{\pt x}\dfr{\pt^3 u}{\pt x^3}\Big] + \mathcal{O}(h^5),\\[2mm]
\beta^{(1)} = h^2 {(\dfr{\pt u}{\pt x})}^2 - h^3\dfr{\pt u}{\pt x}\dfr{\pt^2 u}{\pt x^2} + h^4\Big[\dfr 43 {(\dfr{\pt^2 u}{\pt x^2})}^2 - \dfr 14\dfr{\pt u}{\pt x}\dfr{\pt^3 u}{\pt x^3}\Big] + \mathcal{O}(h^5),\\[2mm]
\beta^{(0)} = h^2 {(\dfr{\pt u}{\pt x})}^2,
\eda
\end{equation}
where $ \pt u/\pt x $, $ \pt^2 u/\pt x^2 $ and $ \pt^3 u/\pt x^3 $ are evaluated at $ x=0 $. 
The linear weights of each stencil are
\begin{equation}
\bga{l}
 d^{(0)} = h^2,\ \ d^{(1)} = h,\ \ d^{(2)} = 1 - d^{(0)} - d^{(1)}.
\eda
\end{equation}
Then, we calculate the weights of the stencils by
\begin{equation}\label{eq:weight}
\al^{(r)} = \dfr{d^{(r)}}{{(\varepsilon+\beta^{(r)})}^2},\ \ \ \
\omega^{(r)} = \dfr{\al^{(r)}}{\sum_{l=0}^2{\al^{(l)}}}.
\end{equation}
Finally, we obtain
\begin{equation}\label{eq:weno-ilw-un1}
\bga{ll}
  \bar{ u}_{-1} = \d\sum_{r=0}^2\omega^{(r)}\bar{ u}_{-1}^{(r)},& \bar{ u}_{-2} = \d\sum_{r=0}^2\omega^{(r)}\bar{ u}_{-2}^{(r)},\\
  \Delta  u_{-1} = \d\sum_{r=0}^2\omega^{(r)}\Delta  u_{-1}^{(r)},& \Delta  u_{-2} = \d\sum_{r=0}^2\omega^{(r)}\Delta  u_{-2}^{(r)}.
\eda
\end{equation}
\vspace{2mm}

\subsection{Inflow boundary condition treatment at intermediate stages}\label{sec:inter-stage}
As pointed out by many researchers  e.g., in \cite{BD-ACC,BD-ACC-2}, the direct use of exact boundary conditions at intermediate stages in the process of multi-stage time discretizations will cause defects of the numerical accuracy. The same argument applies in the current study. We need to treat the boundary condition properly at the intermediate stage.  \\

Our strategy is made as follows.  We first specify to the advancing formula \eqref{eq:1d-al2} for the two-stage fourth-order scheme in the leftmost control volume $I_0$ as
\begin{equation}
\begin{array}{rl}
\bar{u}^{n+1}_0 & =\bar u_0^n -\dfr{ k }{h} \left[f_{\frac 12}^{4th} -f_{-\frac 12}^{4th}\right]\\[3mm]
 &= \bar u_0^n -\dfr{1}{h} \Big\{ k\left[f(u_{\frac 12}^n) -f(u_{-\frac 12}^n)\right]\\[3mm]
 & \ \ \ \ \ \  \ \   +\dfr{ k^2 }{6}\left[f^\prime(u^{n}_{\frac 12})(\dfr{\pt u}{\pt t})^{n}_{\frac 12} -  f^\prime( u^{n}_{-\frac 12})(\dfr{\pt u}{\pt t})^{n}_{-\frac 12}\right]\\[3mm]
&\ \ \ \ \ \ \ \ +\dfr{ k^2 }{3}\left[f^\prime( u^{n+\frac 12}_{\frac 12})(\dfr{\pt u}{\pt t})^{n+\frac 12}_{\frac 12} - f^\prime(u^{n+\frac 12}_{-\frac 12})(\dfr{\pt u}{\pt t})^{n+\frac 12}_{-\frac 12}\right]  \Big\}. 
 \end{array}
 \end{equation}
Using the governing equation \eqref{eq:govern-eq} to replace the temporal derivatives by the corresponding spatial ones, we obtain
\begin{equation}\label{eq:1d-al2-dire}
\begin{array}{rl} 
\bar{u}^{n+1}_0&= \bar u_0^n -\dfr{1}{h} \Big\{ k\left[f(u_{\frac 12}^n) - f(u_{-\frac 12}^n)\right]\\[3mm]
 & \ \ \ \ \ \  \ \ -\dfr{ k^2 }{6}\left[(f^\prime(u^{n}_{\frac 12}))^2(\dfr{\pt u}{\pt x})^{n}_{\frac 12} - (f^\prime(u^{n}_{-\frac 12}))^2(\dfr{\pt u}{\pt x})^{n}_{-\frac 12}\right] \\[3mm]
& \ \ \ \ \  \ \ -\dfr{ k^2 }{3}\left[(f^\prime(u^{n+\frac 12}_{\frac 12}))^2(\dfr{\pt u}{\pt x})^{n+\frac 12}_{\frac 12} - (f^\prime(u^{n+\frac 12}_{-\frac 12}))^2(\dfr{\pt u}{\pt x})^{n+\frac 12}_{-\frac 12}\right] \Big\}.
\end{array}
\end{equation}
The difficulty results from the presence of $ ({\pt u/\pt x})^{n+\frac 12}_{\frac 12} $ and $ ({\pt u/\pt x})^{n+\frac 12}_{-\frac 12} $  evaluated at the intermediate stage $t=t^{n+\frac 12}$. In fact, the interpolation are stated in \eqref{eq:HWENO-general} and specified to this boundary control volume as,
\begin{equation}\label{eq:derv-recon}
\begin{array}{l}
  \left(\dfr{\pt u}{\pt x}\right)^{n+\frac 12}_{-\frac 12} = \dfr{1}{12h}(\bar{u}^{n+\frac 12}_{-2} - 15\bar{u}^{n+\frac 12}_{-1} + 15\bar{u}^{n+\frac 12}_{0} - \bar{u}^{n+\frac 12}_{1}), \\[3mm]
  \left(\dfr{\pt u}{\pt x}\right)^{n+\frac 12}_{\frac 12} = \dfr{1}{12h}(\bar{u}^{n+\frac 12}_{-1} - 15\bar{u}^{n+\frac 12}_{0} + 15\bar{u}^{n+\frac 12}_{1} - \bar{u}^{n+\frac 12}_{2}),
\end{array}
\end{equation}
where $ \bar{u}^{n+\frac 12}_{-1} $ and $ \bar{u}^{n+\frac 12}_{-2} $ are determined in \eqref{eq:ILW-result-u}. Substituting \eqref{eq:ILW-result-u} into \eqref{eq:derv-recon} gives their direct reconstructions
\begin{equation}\label{eq:derv-recon-dire}
\begin{array}{l}
  \left(\dfr{\pt u}{\pt x}\right)^{n+\frac 12}_{-\frac 12} = -(f'(g(t^{n+\frac{1}{2}})))^{-1}g^{\prime}(t^{n+\frac{1}{2}}),\\[3mm]
  \left(\dfr{\pt u}{\pt x}\right)^{n+\frac 12}_{\frac 12} = \dfr{1}{48h}\left[-49\bar{u}^{n+\frac{1}{2}}_{0} + 59\bar{u}^{n+\frac{1}{2}}_{1} - 4\bar{u}^{n+\frac{1}{2}}_{2} \right. \\[3mm]
   \ \ \ \ \ \ \ \ \ \ \ \ \ \ \ \ \ \ \ \ \ \ \ \left. - 6g(t^{n+\frac{1}{2}}) + 6h(f'(g(t^{n+\frac{1}{2}})))^{-1}g^{\prime}(t^{n+\frac{1}{2}})\right].
\end{array}
\end{equation}
The errors in the above formulae come from two sources:  the interpolation error and the truncation error carried by $ \bar{u}^{n+\frac 12}_0 $, $ \bar{u}^{n+\frac 12}_1 $ and $ \bar{u}^{n+\frac 12}_2 $ within  second-order accuracy. The latter one can be analyzed carefully as follows.

Recall the definition of the numerical flux $f^{(1)}_{j+\frac 12}$ in \eqref{eq:1d-al},
\begin{equation}
\begin{array}{rl}
  \dfr k2 f^{(1)}_{j+\frac 12} &= \dfr k2 f(u_{j+\frac 12}^n) + \dfr{ k^2 }{8} (\dfr{\pt f}{\pt t})^{n}_{j+\frac 12}\\[3mm]
& = \d\int^{t^{n+\frac 12}}_{t^n}f(u(x_{j+\frac 12},t))dt - \dfr{ k ^3}{48}\dfr{\pt^2 f}{\pt t^2}(x_{j+\frac 12},t^n) + \mathcal{O}( k ^4).
\end{array}
\end{equation}
Then we find that the intermediate value $ \bar{u}^{n+\frac{1}{2}}_{j} $ in \eqref{eq:derv-recon-dire} bears the truncation error of order $ \mathcal{O}(k^3) $,  compared to the exact solution,  as
\begin{equation}\label{eq:inter-err-0}
\begin{array}{rl}
  \bar{u}^{n+\frac{1}{2}}_{j} & = \bar{u}^{n}_j - \dfr{ k }{2h}\left(f^{(1)}_{j+\frac 12}-f^{(1)}_{j-\frac 12}\right)\\[3mm]
 & = \bar{u}^{n}_j - \dfr{1}{h}\Big[\d\int^{t^{n+\frac 12}}_{t^n}f(u(x_{j+\frac 12},t))dt - \dfr{ k ^3}{48}\dfr{\pt^2 f}{\pt t^2}(x_{j+\frac 12},t^n)\\[3mm]
& \ \ \ \ \ \ \ \ \ - \d\int^{t^{n+\frac 12}}_{t^n}f(u(x_{j-\frac 12},t))dt + \dfr{ k ^3}{48}\dfr{\pt^2 f}{\pt t^2}(x_{j-\frac 12},t^n) + \mathcal{O}( k ^5)\Big]\\[3mm]
 & = \overline{u(\cdot,t^{n+\frac 12})}_{j} + \dfr{ k ^3}{48h}\Big[\dfr{\pt^2 f}{\pt t^2}(x_{j+\frac 12},t^n) - \dfr{\pt^2 f}{\pt t^2}(x_{j-\frac 12},t^n)\Big] + \mathcal{O}( k ^4),
\end{array}
\end{equation}
where the difference in the bracket provides  an $ \mathcal{O}(h) $ term and
\begin{equation*}
\overline{u(\cdot,t^{n+\frac 12})}_{j}=\dfr{1}{h}\displaystyle\int_{I_{j}}u(x,t^{n+\frac{1}{2}})dx, \ \ j=0,1,2.
\end{equation*}
For $j=0,1,2$, we replace the difference quotients in \eqref{eq:inter-err-0} by the corresponding differential quotients at $x=0$ to further get
\begin{equation}\label{eq:inter-err-1}
\begin{array}{rl}
\bar{u}^{n+\frac{1}{2}}_{j} & = \overline{u(\cdot,t^{n+\frac 12})}_{j} + \dfr{ k ^3}{48}\dfr{\pt^3 f}{\pt t^2 \pt x}(0,t^n) + \mathcal{O}( k ^4)\\[3mm]
& = \overline{u(\cdot,t^{n+\frac 12})}_{j} - \dfr{ k ^3}{48}\dfr{\pt^3 u}{\pt t^3}(0,t^n) + \mathcal{O}( k ^4)\\[3mm]
& = \overline{u(\cdot,t^{n+\frac 12})}_{j} - \dfr{ k ^3}{48}g^{\prime\prime\prime}(t^n) + \mathcal{O}( k ^4)\\[3mm]
& = \overline{u(\cdot,t^{n+\frac 12})}_{j} - \dfr{ k ^3}{48}g^{\prime\prime\prime}(t^{n+\frac 12}) + \mathcal{O}( k ^4),
\end{array}
\end{equation}
where the second equality is derived by replacing the spatial derivative by the corresponding temporal derivative through the governing equation in \eqref{eq:IBVP}. Substituting \eqref{eq:inter-err-1} into \eqref{eq:derv-recon-dire} gives us
\begin{equation}\label{eq:derv-recon-dire-error}
\begin{array}{l}
  \left(\dfr{\pt u}{\pt x}\right)^{n+\frac 12}_{-\frac 12} = -(f'(g(t^{n+\frac{1}{2}})))^{-1}g^{\prime}(t^{n+\frac{1}{2}})=\dfr{\pt u}{\pt x}(x_{-\frac 12},t^{n+\frac 12}),\\[3mm]
  \left(\dfr{\pt u}{\pt x}\right)^{n+\frac 12}_{\frac 12} = \dfr{1}{48h}\Big[-49\overline{u(\cdot,t^{n+\frac 12})}_{0} + 59\overline{u(\cdot,t^{n+\frac 12})}_{1} - 4\overline{u(\cdot,t^{n+\frac 12})}_{2}\\[3mm]
\ \ \ \ \ \ \ \ \ \ \ \ \ \ \ \ \ \ \ \ \ \ \ \ \ \ - 6g(t^{n+\frac{1}{2}}) + 6h(f'(g(t^{n+\frac{1}{2}})))^{-1}g^{\prime}(t^{n+\frac{1}{2}})\\[3mm]
\ \ \ \ \ \ \ \ \ \ \ \ \ \ \ \ \ \ \ \ \ \ \ \ \ \ - \dfr{6k^3}{48}g^{\prime\prime\prime}(t^{n+\frac 12})\Big],\\[3mm]
\ \ \ \ \ \ \ \ \ \ \ \ \ \ \ =\dfr{\pt u}{\pt x}(x_{\frac 12},t^{n+\frac 12}) + \mathcal{O}(h^3) - \dfr{1}{48h}\dfr{6k^3}{48}g^{\prime\prime\prime}(t^{n+\frac 12}).
\end{array}
\end{equation}
The error $\mathcal{O}(h^3)$ in the second equation in \eqref{eq:derv-recon-dire-error} comes from the interpolation approximation while $\dfr{1}{48h}\dfr{6k^3}{48}g^{\prime\prime\prime}(t^{n+\frac 12})$ is induced from \eqref{eq:inter-err-1}.  It is easy to see that substituting \eqref{eq:derv-recon-dire-error} into \eqref{eq:1d-al2-dire} leads to
\begin{equation*}
  \bar{u}_0^{n+1}=\dfr{1}{h}\displaystyle\int_{I_{0}}u(x,t^{n+1})dx + \mathcal{O}(k^3),
\end{equation*}
which means that the numerical scheme is only third-order accurate if the exact value of the boundary data is input directly in $I_0$. 

In order to restore the fourth-order accuracy of the two-stage fourth-order scheme, we need to eliminate the $\mathcal{O}(k^2)$ error in \eqref{eq:derv-recon-dire-error}. For this purpose,  we use
\begin{equation}\label{eq:derv-recon-dire-m}
\begin{array}{l}
  \left(\dfr{\pt u}{\pt x}\right)^{n+\frac 12}_{-\frac 12} = -(f^\prime(g(t^{n+\frac{1}{2}})))^{-1}(g^{\prime})^{n+\frac{1}{2}},\\[3mm]
  \left(\dfr{\pt u}{\pt x}\right)^{n+\frac 12}_{\frac 12} = \dfr{1}{48h}\left[-49\bar{u}^{n+\frac{1}{2}}_{0} + 59\bar{u}^{n+\frac{1}{2}}_{1} - 4\bar{u}^{n+\frac{1}{2}}_{2} \right. \\[3mm]
   \ \ \ \ \ \ \ \ \ \ \ \ \ \ \ \ \ \ \ \ \ \ \ \left. - 6g^{n+\frac{1}{2}} + 6h(f^\prime(g(t^{n+\frac{1}{2}})))^{-1}(g^{\prime})^{n+\frac{1}{2}}\right],
\end{array}
\end{equation}
to reconstruct $({\pt u}/{\pt x})^{n+\frac 12}_{\pm\frac 12}$ instead of \eqref{eq:derv-recon-dire}. Here the exact boundary values $g(t^{n+\frac 12})$ and $g^\prime(t^{n+\frac 12})$ used in \eqref{eq:derv-recon-dire} are replaced by
\begin{equation}\label{eq:inter-err-g}
  \begin{array}{l}
g^{n+\frac 12} = g(t^{n+\frac 12}) + e^{(0)} k^3,\\[3mm]
{(g^{\prime})}^{n+\frac 12} = g^{\prime}(t^{n+\frac 12}) + e^{(1)} k^2.\\
  \end{array}
\end{equation}
The terms $ e^{(0)} k^3 $ and $ e^{(1)} k^2 $ are introduced into $ g^{n+\frac 12} $ and $ {(g^{\prime})}^{n+\frac 12} $ to cancel the truncation errors carried by $ \bar{u}^{n+\frac 12}_0 $, $ \bar{u}^{n+\frac 12}_1 $ and $ \bar{u}^{n+\frac 12}_2 $ and they are to be determined in the following. By combining \eqref{eq:inter-err-1}, \eqref{eq:derv-recon-dire-m} and \eqref{eq:inter-err-g}, we have
\begin{equation}\label{eq:derv-recon-error-m-1}
\begin{array}{rl}
  \left(\dfr{\pt u}{\pt x}\right)^{n+\frac 12}_{-\frac 12}& = -(f^\prime(g(t^{n+\frac 12})))^{-1}{(g^{\prime})}^{n+\frac 12}\\[3mm]
&= -(f^\prime(g(t^{n+\frac 12})))^{-1}g^{\prime}(t^{n+\frac 12}) - (f^\prime(g(t^{n+\frac 12})))^{-1} e^{(1)} k^2\\[3mm]
& =\dfr{\pt u}{\pt x}(x_{-\frac 12},t^{n+\frac 12}) -(f^\prime(g(t^{n+\frac 12})))^{-1} e^{(1)} k^2,
\end{array}
\end{equation}
and
\begin{equation}
\begin{array}{rl}
  \left(\dfr{\pt u}{\pt x}\right)^{n+\frac 12}_{\frac 12}&= \dfr{1}{48h}\left[-49\bar{u}^{n+\frac 12}_{0} + 59\bar{u}^{n+\frac 12}_{1} - 4\bar{u}^{n+\frac 12}_{2} \right.\\[3mm]
&\  \ \ \ \  \ \ \ \ \ \ \  \left. - 6g^{n+\frac 12} + 6h(f^\prime(g(t^{n+\frac 12})))^{-1}{(g^{\prime})}^{n+\frac 12}\right]\\[3mm]
 &= \dfr{1}{48h}\left[-49\overline{u(\cdot,t^{n+\frac 12})}_{0} + 59\overline{u(\cdot,t^{n+\frac 12})}_{1} - 4\overline{u(\cdot,t^{n+\frac 12})}_{2} \right.\\[3mm]
&\ \ \ \ \ \   \ \ \ \ \ \ -6g(t^{n+\frac 12}) +6h(f^\prime(g(t^{n+\frac 12})))^{-1}g^{\prime}(t^{n+\frac 12})\\[3mm]
&\ \ \ \ \ \ \ \ \  \ \ \left.  -\dfr{6}{48}k^3g^{\prime\prime\prime}(t^{n+\frac 12}) - 6e^{(0)}k^3 + 6(f^\prime(g(t^{n+\frac 12})))^{-1}e^{(1)}hk^2\right]. 
\end{array}
\end{equation}
This can be evaluated as 
\begin{equation}\label{eq:derv-recon-error-m-2}
\begin{array}{rl}
\left(\dfr{\pt u}{\pt x}\right)^{n+\frac 12}_{\frac 12}&=\dfr{\pt u}{\pt x}(x_{\frac 12},t^{n+\frac 12}) + \mathcal{O}(h^3)\\[3mm]
&\ \ \ \ \  + \dfr{k^3}{48h}\left[-\dfr{6}{48}g^{\prime\prime\prime}(t^{n+\frac 12}) -6e^{(0)}+6(f^\prime(g(t^{n+\frac 12})))^{-1}e^{(1)}\dfr hk\right].\\[3mm]
\end{array}
\end{equation}
In order to eliminate the $\mathcal{O}(k^2)$ terms in both \eqref{eq:derv-recon-error-m-1} and \eqref{eq:derv-recon-error-m-2}, we require
\begin{equation}\label{eq:inter-err-equations}
\begin{array}{r}
 f^\prime(g(t^{n+\frac{1}{2}}))^{-1}e^{(1)} = 0,\\[2mm]
  -\dfr{6}{48}g^{\prime\prime\prime}(t^{n+\frac 12}) -6e^{(0)}+6(f^\prime(g(t^{n+\frac 12})))^{-1}e^{(1)}\dfr hk = 0.
\end{array}
\end{equation}
Solving the above system of linear equations with respect to $ e^{(0)} $ and $ e^{(1)} $ provides
\begin{equation}\label{eq:inter-err-linear-results}
\begin{array}{l}
  e^{(0)} = -\dfr{1}{48}g^{\prime\prime\prime}(t^{n+\frac 12}),\\[2mm]
  e^{(1)} = 0.
\end{array}
\end{equation}
Therefore at the intermediate stage $ t^{n+\frac 12} $, the boundary values used in numerical computations are modified as
\begin{equation}\label{eq:inter-boundary-value}
\begin{array}{l}
  g^{n+\frac{1}{2}} = g(t^{n+\frac{1}{2}}) + e^{(0)}  k ^3 = g(t^{n+\frac{1}{2}}) -\dfr{ k ^3}{48}g^{\prime\prime\prime}(t^{n+\frac 12}),\\[2mm]
  {(g^{\prime})}^{n+\frac{1}{2}} = g^\prime(t^{n+\frac{1}{2}}).\\
\end{array}
\end{equation}
\vspace{2mm}

\begin{rem}
The final intermediate boundary values \eqref{eq:inter-boundary-value} are consistent with those  given in \cite{BD-ACC-2} where the author indicates that the boundary values used at the intermediate stage should have the same truncation errors as those in the interior cells when multi-stage methods are used for initial and boundary value hyperbolic problems.
\end{rem}
\vspace{0.2cm}

\subsection{Outflow boundary condition treatment}\label{sec:outflow}
For the initial and boundary problem  \eqref{eq:IBVP}, we set  $ x_{M+\frac 12}=1 $ as an outflow boundary, at which  no boundary condition is required theoretically. However, we have to interpolate it to obtain the required values $ \bar{u}_{M+1} $, $ \bar{u}_{M+2} $, $ \Delta u_{M+1} $ and $ \Delta u_{M+2} $ in ghost cells in order to implement the scheme.
Since the signal propagates out of the computational domain through the boundary $ x=1 $, the extrapolation can be  used to construct the data in the  ghost cells  $I_{M+1}$ and $I_{,M+2}$. Here the extrapolation developed in \cite{ILW} is briefly summarized.

 In order to achieve the fourth-order accuracy, a cubic polynomial is constructed by using $\bar{u}_{M-3}$, $ \bar{ u}_{M-2} $, $ \bar{ u}_{M-1} $ and $ \bar{ u}_{M} $, which is
\begin{equation}\label{eq:extra-4-poly}
\bga{l}
  q(x) = \dfr{\bar{u}_{M}-3\bar{u}_{M-1}+3\bar{u}_{M-2}-\bar{ u}_{M-3}}{6h^3}(x-1)^3\\[3mm]
\ \ \ \ \ \ \ \ + \dfr{5\bar{u}_{M}-13\bar{u}_{M-1}+11\bar{u}_{M-2}-3\bar{u}_{M-3}}{4h^2}(x-1)^2\\[3mm]
\ \ \ \ \ \ \ \  + \dfr{35\bar{u}_{M}-69\bar{u}_{M-1}+45\bar{u}_{M-2}-11\bar{u}_{M-3}}{12h}(x-1)\\[3mm]
\ \ \ \ \ \ \ \  + \dfr{25\bar{u}_{M}-23\bar{u}_{M-1}+13\bar{u}_{M-2}-3\bar{u}_{M-3}}{12}.
\eda
\end{equation}
This gives the values
\begin{equation}\label{eq:extra-4-res}
\bga{l}
  \bar{u}_{M+1} = 4\bar{u}_{M} - 6\bar{u}_{M-1} + 4\bar{u}_{M-2} - \bar{u}_{M-3},\\[3mm]
  \bar{u}_{M+2} = 10\bar{u}_{M} - 20\bar{u}_{M-1} + 15\bar{u}_{M-2} - 4\bar{u}_{M-3},\\[3mm]
  \Delta u_{M+1} = \dfr{26\bar{u}_{M} - 57\bar{u}_{M-1} + 42\bar{u}_{M-2} - 11\bar{u}_{M-3}}{6},\\[3mm]
  \Delta u_{M+2} = \dfr{47\bar{u}_{M} - 114\bar{u}_{M-1} + 93\bar{u}_{M-2} - 26\bar{u}_{M-3}}{6}.
\eda
\end{equation}
If there is a discontinuity in either $ I_{M-3} $, $ I_{M-2} $, $ I_{M-1} $ or $ I_{M} $, a WENO-type stencil selection can be applied. Details can be found in \cite{high-order-ILW}.
\vspace{0.2cm}

\section{Numerical boundary condition treatment for hyperbolic systems}\label{sec:system}
This section turns to investigate numerical treatment for boundary conditions of hyperbolic systems. We focus on hyperbolic conservation laws. As far as source terms are included, the corresponding extensions are automatically made.  Systems of  hyperbolic conservation laws in one dimension read
\begin{equation}\label{eq:system}
  \dfr{\pt\bu}{\pt t} + \dfr{\pt\bbf(\bu)}{\pt x} = 0,
\end{equation}
where $ \bu $ has $ m $ components $ \bu=( u_1, \dots, u_m )^\top$, $ x\in (0,1) $ and $ t > 0 $, $\bbf(\bu)$ is the associated flux vector function. In this section, all bold letters represent vectors. In contrast with the scalar case in Section \ref{sec:scalar},  the case of systems is more complicated since inflow and outflow signals are coupled. The basic idea of numerical treatments originates  from \cite{EB-consv-FD,ILW}.  Inflow signals and  outflow signals are separated  to be treated using the characteristic decomposition so that  the technique of inflow and outflow boundary condition treatments in  Section \ref{sec:scalar} can be applied, respectively.

Denote by $ q^{(k)}:= \dfr{\pt^k q}{\pt x^k} $ the $ k $-th order spatial derivative of a certain quantity $ q $ at the boundary, $k=0,1$, and denote by  $\dfr{\pt\bbf(\bu)}{\pt\bu}$ the Jacobian matrix of $\bbf(\bu)$. Assume that there are $ m $ real eigenvalues of $\dfr{\pt\bbf}{\pt\bu}(\bu^{(0)})$
\begin{equation}\label{eq:eigenvalue}
  \lambda_1 \leq \dots \leq \lambda_r \leq 0 < \lambda_{r+1} \leq \dots \leq \lambda_{m},
\end{equation}
 and a   complete set of corresponding left eigenvectors $\mathbf{L}_i(\bu^{(0)})$,
\begin{equation}
 \bL_i(\bu^{(0)})^\top \dfr{\pt\bbf}{\pt\bu}(\bu^{(0)}) =\la_i\bL_i(\bu^{(0)})^\top.
\end{equation}
Thus there are $ m-r $ out-going signals from the boundary $ x=0 $ and the number of boundary conditions should  equal to the number of out-going characteristics, symbolically denoted as
\begin{equation}\label{eq:boundary-condition}
  B_i(\bu(0,t)) = \psi_i(t),  \ \ \ \ \ i = r+1,\dots,m.
\end{equation}
We refer to \cite{ibvp} for details. The similar boundary conditions can be imposed when the right boundary is considered.

In the following part of this section, we first propose the general procedure for one-dimensional conservation systems and then specify to several  examples, including the solid-wall boundary condition for the one-dimensional Euler equations,   the subsonic inflow and the subsonic outflow boundary conditions involved in computations of the nuzzle flow problem.  In the last subsection, we extend our boundary condition treatment to the two-dimensional case by implementing it for the solid-wall boundary condition.\\

\subsection{General framework for hyperbolic systems}\label{subsec:system-general}
At first, introduce  characteristic variables
\begin{equation}\label{eq:decomp}
  (\bar{w}_i)_j = \bL_i \cdot \bar{\bu}_j,\ i = 1,\dots,r,
\end{equation}
where $ \bar{\bu}_j $ is the cell average of $ \bu $ over the cell $ I_j $. Theoretically, $\bL_i(\bu^{(0)})$ should be used here.  However in practice, we can use $\bL_i(\bar{\bu}_0)$ in the expression  \eqref{eq:decomp} if $\bL_i(\bu^{(0)})$ is not available.  Since $ w_1,\dots,w_r $ are out-going signals, the extrapolation technique in Subsection \ref{sec:outflow} is applied to obtain $ w_i^{(0)} $ and $ w_i^{(1)} $ for $ i=1,\dots,r $. For example, as the solution is smooth close to the boundary $ x=0 $, similar to \eqref{eq:extra-4-poly}, an interpolation polynomial with respect to $w_i$ near $x=0$ can be constructed which leads to
\begin{equation}\label{eq:extra-bod-res}
\bga{l}
  w_i^{(0)} = \dfr{1}{12}\Big[-25 (\bar{w}_i)_0 + 23 (\bar{w}_i)_1 - 13 (\bar{w}_i)_2 + 3 (\bar{w}_i)_3\Big],\\[2mm]
  w_i^{(1)} = \dfr{1}{12h}\Big[-35 (\bar{w}_i)_0 + 69 (\bar{w}_i)_1 - 45 (\bar{w}_i)_2 + 11 (\bar{w}_i)_3\Big].
\eda
\end{equation}
Then we solve the system
\begin{equation}\label{eq:general-0}
\left\{
\bga{l}
\bL_1\cdot\bu^{(0)} = w^{(0)}_1,\\
\ \ \ \ \ \vdots\\
\bL_r\cdot\bu^{(0)} = w^{(0)}_r, \\
B_{r+1}(\bu^{(0)}) = \td\psi_{r+1},\\
\ \ \ \ \ \vdots\\
B_m(\bu^{(0)}) = \td\psi_m,
\eda
\right.
\end{equation}
to obtain $\bu^{(0)}:=\bu(0,t)$. Here $ w_i^{(0)} $ is  obtained by \eqref{eq:extra-bod-res}. For $ i\in\{r+1,\dots,m\} $, $\td\psi_i$ is the modified boundary data defined by $\td\psi_i=\psi_{i}(t^{n+\frac 12})-\dfr{k^3}{48}\psi_{i}^{\prime\prime\prime}(t^{n+\frac 12})$ at the intermediate stage $t=t^{n+\frac 12} $ and $\td\psi_i=\psi_{i}(t^n)$ at $t= t^n$.

In order to derive $\bu^{(1)} := {\pt\bu}/{\pt x}(0,t) $, we carry out the following manipulation.  We take temporal derivatives on both sides of the boundary conditions in \eqref{eq:boundary-condition} to  yield
\begin{equation}\label{eq:general-dt}
  \nabla_{\bu}B_i(\bu^{(0)})^\top \dfr{\pt \bu}{\pt t}(0,t) = \psi_i^\prime(t), \ \ \ \ \ \  i = r+1,\dots,m.
\end{equation}
Then the governing equations \eqref{eq:system} are adopted to convert temporal derivatives to spatial ones as
\begin{equation}\label{eq:general-dx}
\nabla_{\bu}B_i(\bu^{(0)})^\top \Big[- \dfr{\pt\bbf}{\pt\bu}(\bu^{(0)}) \ \dfr{\pt \bu}{\pt x}(0,t) \Big]= \psi_i^\prime(t), \ \ \ \  i = r+1,\dots,m.
\end{equation}
Furthermore, we have
\begin{equation}
  \bL_i\cdot\bu^{(1)} = w^{(1)}_i, \ \ \ \  \ i=1,\dots,r,
\end{equation}
by taking spatial derivatives on both sides of the first $ r $ equations in \eqref{eq:general-0}. Thus  we derive the following linear system
\begin{equation}\label{eq:general-1}
\left\{
\bga{l}
\bL_1\cdot\bu^{(1)} = w^{(1)}_1,\\
\ \ \ \ \ \ \ \ \ \ \ \vdots \\
\bL_r\cdot\bu^{(1)} = w^{(1)}_r,\\
\nabla_{\bu}B_{r+1}(\bu^{(0)})^\top \ \dfr{\pt\bbf}{\pt\bu}(\bu^{(0)}) \ \bu^{(1)} = -\psi_{r+1}^\prime,\\
\ \ \ \ \ \ \ \ \ \ \ \vdots \\
\nabla_{\bu}B_m(\bu^{(0)})^\top \ \dfr{\pt\bbf}{\pt\bu}(\bu^{(0)}) \ \bu^{(1)} = -\psi_m^\prime,
\eda
\right.
\end{equation}
with $ \bu^{(1)} $ as the unknown and $ w_i^{(1)} $ defined in \eqref{eq:extra-bod-res}.

After $ \bu^{(0)}=[u^{(0)}_1,\dots,u^{(0)}_m]^\top $ and $ \bu^{(1)}=[u^{(1)}_1,\dots,u^{(1)}_m]^\top $ are obtained by solving linear systems \eqref{eq:general-0} and \eqref{eq:general-1}, for each pair $ u^{(0)}_i $ and $ u^{(1)}_i $, $ i=1,\dots,m $, we construct polynomials $ p_i^{(r)} $ on stencils $ S_i^{(r)} $, $ r=0,1,2 $, under the conditions $ p_i^{(r)}(0)=u^{(0)}_i $ and $ {[p_i^{(r)}]}^\prime(0)=u^{(1)}_i $, as in Subsection \ref{sec:inflow}. Then $ (\bar{u_i})_{-1} $, $ (\bar{u_i})_{-2} $, $ (\Delta u_i)_{-1} $ and $ (\Delta u_i)_{-2} $ are defined as in \eqref{eq:weno-ilw-un1}. Finally, we obtain
\begin{equation}
\bga{cc}
  \bar{\bu}_{-1}=\left[
    \begin{array}{c}
      (\bar{u_1})_{-1}\\
       \vdots\\
      (\bar{u_m})_{-1}
    \end{array}
  \right], &
  \bar{\bu}_{-2}=\left[
    \begin{array}{c}
      (\bar{u_1})_{-2}\\
       \vdots\\
      (\bar{u_m})_{-2}
    \end{array}
  \right], \\[4mm]
  \Delta{\bu}_{-1}=\left[
    \begin{array}{c}
      (\Delta{u_1})_{-1}\\
       \vdots\\
      (\Delta{u_m})_{-1}
    \end{array}
  \right], &
  \Delta{\bu}_{-2}=\left[
    \begin{array}{c}
      (\Delta{u_1})_{-2}\\
       \vdots\\
      (\Delta{u_m})_{-2}
    \end{array}
  \right]
\eda
\end{equation}
for practical computations.\\

\vspace{2mm}

\subsection{Solid-wall boundary condition for the one-dimensional Euler equations}\label{subsec:wall}
In this subsection, we will practically apply   \eqref{eq:general-0} and \eqref{eq:general-1} for the solid-wall boundary condition of the one-dimensional compressible  Euler equations,
\begin{equation}\label{eq:Euler}
 \bu = (\rho,  \rho v, \rho E)^\top, \ \ \ \
  \bbf(\bu) =( \rho v, \rho v^2+p,  v(\rho E + p))^\top,
\end{equation}
where $\rho, v, p$ are the density, velocity and pressure, $ E = \dfr 12 v^2 +e $ is the total energy with the internal energy $ e = e(p,\rho) $, which is given  in terms of the equation of state (EOS). In this paper, we just consider the case of polytropic gases with $e = \dfr{p}{(\gamma-1)\rho}$,  $\gm>1$.
Then the Jacobian matrix of the flux function in terms of the conservative variables $\bu$ takes
\begin{equation}\label{eq:jacobi-Euler}
  \dfr{\pt \bbf}{\pt \bu} = \left[
    \begin{array}{ccc}
      0&1&0\\
\dfr{\gamma-3}{2}v^2 & (3-\gamma)v & \gamma-1\\
\dfr{\gamma-2}{2}v^3-\dfr{1}{\gamma-1}vc^2&\dfr{3-2\gamma}{2}v^2+\dfr{1}{\gamma-1}c^2&\gamma v
    \end{array}
  \right].
\end{equation}
It has three eigenvalues  $ \lambda_1=v-c $, $ \lambda_2=v $ and $ \lambda_3=v+c $, where $ c = \sqrt{\gm p/\rho} $ is the sound speed. The corresponding left eigenvectors are
\begin{equation}\label{eq:lefteigen-Euler}
  \bga{l}
\bL_1 = \dfr{1}{2c^2}\left[\dfr{\gamma-1}{2}v^2 + vc, -(\gamma-1)v-c, \gamma-1\right]^\top,\\[3mm]
\bL_2 =-\dfr{1}{c^2}\left[\dfr{\gamma-1}{2}v^2 - c^2, -(\gamma-1)v, \gamma-1\right]^\top,\\[3mm]
\bL_3 = \dfr{1}{2c^2}\left[\dfr{\gamma-1}{2}v^2 - vc, -(\gamma-1)v+c, \gamma-1\right]^\top.
  \eda
\end{equation}

Assume that $x=0$ is a  solid wall, on which the flow velocity is zero, $ v(0,t) = 0 $ which leads to
\begin{equation}\label{eq:wall}
  u_2(0,t) = \rho(0,t) v(0,t) = 0.
\end{equation}
So the system \eqref{eq:general-0}, when specified to the current case,  is
\begin{equation}\label{eq:wall-0}
\left\{
\bga{l}
\bL_1\cdot\bu^{(0)} = w^{(0)}_1,\\[2mm]
\bL_2\cdot\bu^{(0)} = w^{(0)}_2,\\[2mm]
w^{(0)}_2 = 0.
\eda
\right.
\end{equation}
The modification term in \eqref{eq:general-0} is zero for the present case since $ \psi $ takes a constant value at the boundary and therefore $ \psi^{\prime\prime\prime} $ is always zero.

Taking temporal derivatives on both sides of $ u^{(0)}_2 = 0 $ and converting temporal derivatives into spatial ones, we obtain
\begin{equation}
\bga{l}
  0=[0, \ 1, \ 0]\ \dfr{\pt\bbf}{\pt\bu}(\bu^{(0)})\ \bu^{(1)}\\[3mm]
\ \ = [\dfr{\gamma-3}{2}(v^{(0)})^2, \ (3-\gamma)v^{(0)}, \ \gamma-1]\ \bu^{(1)}\\[3mm]
\ \ = [0, \ \ 0, \ \ \gamma-1]\ \bu^{(1)}.
\eda
\end{equation}
The last equality results from the boundary condition  $ v^{(0)} $ = 0. So \eqref{eq:general-1} becomes,  for this specific case,
\begin{equation}\label{eq:wall-1}
\left\{
\bga{l}
\bL_1\cdot\bu^{(1)} = w^{(1)}_1,\\[2mm]
\bL_2\cdot\bu^{(1)} = w^{(1)}_2,\\[2mm]
\left[0, \ \ 0, \ \ \gamma -1 \right] \ \bu^{(1)} = 0.
\eda
\right.
\end{equation}

The numerical example with the solid-wall boundary condition is presented  in  Example 5 when dealing with the Woodward--Colella problem.

\vspace{2mm}

\subsection{Inflow and outflow boundary conditions for the nozzle flow}\label{subsec:in-and-out}
The nozzle flow is ubiquitous  in gas dynamics. See \cite{Ben-Artzi-01} and references therein. The nozzle flow is an IBVP for the Euler equations with the geometry effect resulting from the shape of the duct. The governing equations takes the form
\begin{equation}\label{eq:nozzle}
\dfr{\pt}{\pt t}\left[
\begin{array}{c}
A\rho\\
A\rho v\\
A\rho E
\end{array}
\right]
+ \dfr{\pt}{\pt x}\left[
\begin{array}{c}
A\rho v\\
A(\rho v^2+p)\\
Av(\rho E + p)
\end{array}
\right] = \left[
\begin{array}{c}
0\\
A^\prime p\\
0
\end{array}
\right],
\end{equation}
where $A=A(x)$ is the sectional area of the duct.  In Example 6 of the next section, the duct occupies the computational domain $x\in (0,1)$. The fluid flows into the duct at $ x=0 $ and flows out of the duct at $ x=1 $. Therefore three kinds of boundary conditions are involved in the computations and they  are the subsonic inflow boundary condition, the subsonic outflow boundary condition and the supersonic boundary condition.  As for the supersonic out-going flow at $x=1$, we can simply extrapolate $ \bu $ by using the WENO-type extrapolation component-wise.\\

{\bf (I) Subsonic inflow boundary condition.}  As $|v|<c$ at $x=0$, the inflow is subsonic. That is,   $ v-c < 0 $, $ v+c>v>0 $ at $ x=0 $. Then two boundary conditions are required at this end. In our computation, the inflow pressure and the inflow density are given as
\begin{equation}\label{eq:nozzle-in}
\bga{r}
  (\gamma-1)\Big[u_3(0,t) - \dfr 12 \dfr{{u_2(0,t)}^2}{u_1(0,t)}\Big] = A(0) p_\text{in}(t),\\[2mm]
  u_1(0,t) = A(0) \rho_\text{in}(t).
\eda
\end{equation}
These two equations and  $ \bL_1\cdot\bu^{(0)} = w^{(0)}_1 $ together meet the  form of \eqref{eq:general-0} specifically
\begin{equation}\label{eq:nozzle-in-0}
\left\{
\bga{l}
  \bL_1\cdot\bu^{(0)} = w^{(0)}_1,\\[2mm]
  (\gamma-1)\Big[u_3^{(0)} - \dfr 12 \dfr{{u_2^{(0)}}^2}{u_1^{(0)}}\Big] = A(0) \td p_\text{in},\\[2mm]
  u_1^{(0)} = A(0) \td\rho_\text{in},
\eda
\right.
\end{equation}
where the notations are 
\begin{equation*}
  \td\rho_\text{in} = \rho_\text{in}(t^{n+\frac 12})-\dfr{k^3}{48}\rho_\text{in}^{\prime\prime\prime}(t^{n+\frac 12}), \ \td p_\text{in}=p_\text{in}(t^{n+\frac 12})-\dfr{k^3}{48}p_\text{in}^{\prime\prime\prime}(t^{n+\frac 12}),
\end{equation*}
at the intermediate stage $t=t^{n+\frac 12}$ and
\begin{equation*}
  \td\rho_\text{in} = \rho_\text{in}(t^{n}), \ \td p_\text{in}=p_\text{in}(t^n),
\end{equation*}
at $t=t^n$.

In order to get  $\bu^{(1)}$, we take temporal derivatives on both sides of \eqref{eq:nozzle-in} and  obtain
\begin{equation}\label{eq:nozzle-in-1}
\bga{l}
\ \ \ \Big[\dfr 12 (v^{(0)})^3 - \dfr{1}{\gamma-1}(c^{(0)})^2v^{(0)}, \ -(v^{(0)})^2+\dfr{1}{\gamma-1}(c^{(0)})^2, \ v^{(0)}\Big] \bu^{(1)}\\[3mm]
\ \ \ \ \ \ \ \ \ \ \ \ \ = -\dfr{A(0) p_{\text{in}}^{\prime}}{\gm-1} - A^\prime(0) p^{(0)}v^{(0)},\\[3mm]
\ \ \ [0,\ 1,\ 0] \ \bu^{(1)} = - A(0) \rho_{\text{in}}^{\prime}.
\eda
\end{equation}
The combination of  them and $ \bL_1\cdot\bu^{(1)} = w^{(1)}_1 $ possesses the structure of \eqref{eq:general-1}.\\

{\bf (II) Subsonic outflow boundary condition.}  At the exit of the nozzle $ x=1 $, the out-going flow is subsonic if  $ v-c < 0 $ and $ v+c > v > 0 $.   This means that only the  characteristics associated with $v-c$  are impinging  onto the exit $x=1$  from the exterior  of the computational domain. Therefore just one boundary condition is prescribed theoretically at this end. In the conventional treatment, the outflow pressure at the exit is given, denoted by $ p_{\text{ex}} $. So the boundary condition  is
\begin{equation}\label{eq:nozzle-ex}
  (\gamma-1)\Big[u_3(1,t) - \dfr 12 \dfr{{u_2(1,t)}^2}{u_1(1,t)}\Big] = A(1) p_\text{ex}(t).
\end{equation}
The above equation together with $ \bL_2\cdot\bu^{(0)} = w^{(0)}_2 $ and $ \bL_3\cdot\bu^{(0)} = w^{(0)}_3 $ satisfy the form of \eqref{eq:general-0} as
\begin{equation}\label{eq:nozzle-ex-0}
\left\{
\bga{l}
  (\gamma-1)\Big[u_3^{(0)} - \dfr 12 \dfr{{u_2^{(0)}}^2}{u_1^{(0)}}\Big] = A(1) \td p_\text{ex},\\[2mm]
  \bL_2\cdot\bu^{(0)} = w^{(0)}_2,\\[2mm]
  \bL_3\cdot\bu^{(0)} = w^{(0)}_3,
\eda
\right.
\end{equation}
where $\td p_\text{ex}=p_\text{ex}(t^{n+\frac 12})-\dfr{k^3}{48}p_\text{ex}^{\prime\prime\prime}(t^{n+\frac 12})$ at the intermediate stage $t=t^{n+\frac 12}$ and $\td p_\text{ex}=p_\text{ex}(t^n)$ at $t=t^n$.

Once again, for the linear system of $ \bu^{(1)} $, take temporal derivatives on both sides of \eqref{eq:nozzle-ex} and convert temporal derivatives to spatial ones to get
\begin{equation}\label{eq:nozzle-ex-sub-1}
\bga{l}
\Big[\dfr 12 (v^{(0)})^3 - \dfr{1}{\gamma-1}(c^{(0)})^2v^{(0)}, \ -(v^{(0)})^2+\dfr{1}{\gamma-1}(c^{(0)})^2, \ v^{(0)}\Big] \bu^{(1)}\\[3mm]
\ \ \ \ \ \ \ \ \ \ \ \ \ \ \ \ \ \ \ \ \ \ \ \ \ \ \ \ \ \ \ \ \ \ \ \ \ \ \ \ \ \ \ \ \ \ \ \ \ \ \ \ = -\dfr{A(1) p_{\text{ex}}^{\prime}}{\gm-1} - A^\prime(1) p^{(0)}v^{(0)}.
\eda
\end{equation}
Then, we combine it with $ \bL_2\cdot\bu^{(1)} = w^{(1)}_2 $  and $ \bL_3\cdot\bu^{(1)} = w^{(1)}_3 $  meet the  form of \eqref{eq:general-1}.\\

\subsection{Solid-wall boundary condition for the two-dimensional Euler equations}\label{subsec:wall-2d}
In this subsection, we show how this boundary condition treatment deals with the solid-wall boundary condition of the two-dimensional Euler equations
\begin{equation}\label{eq:Euler-2d}
\dfr{\pt \bu}{\pt t} + \dfr{\pt \bbf(\bu)}{\pt x} + \dfr{\pt \bg(\bu)}{\pt y} = 0
\end{equation}
with
\begin{equation}\label{eq:2D}
 \bu = \left[
    \begin{array}{c}
      \rho\\
      \rho v^x\\
      \rho v^y\\
      \rho E
    \end{array}
  \right] \hspace{8mm}
  \bbf(\bu) = \left[
    \begin{array}{c}
      \rho v^x\\
      \rho (v^x)^2+p\\
      \rho v^xv^y\\
       v^x(\rho E + p)
    \end{array}
  \right] \hspace{8mm}
  \bg(\bu) = \left[
    \begin{array}{c}
      \rho v^y\\
      \rho v^xv^y\\
      \rho (v^y)^2+p\\
       v^y(\rho E + p)
    \end{array}
  \right],
\end{equation}
where $\rho, (v^x,v^y), p$ are the density, velocity and pressure, $ E = \dfr 12 ((v^x)^2+(v^y)^2) +e $ is the total energy with the internal energy  $ e = \dfr{1}{\gamma-1}\dfr p\rho $ for  polytropic gases. Consider the computational domain $ \Omega=\{(x,y) : x > 0 , y\in(y_\text{min},y_\text{max})\} $ with a solid wall $ \Gm=\{(x,y) : x=0, y\in(y_\text{min},y_\text{max})\} $. Here $\bbf$ is the flux normal to $\Gm$ and its Jacobian matrix is
\begin{equation}\label{eq:jacobi-Euler-f}
  \dfr{\pt \bbf}{\pt \bu} = \left[
    \begin{array}{cccc}
      0&1&0&0\\
      (\gm-1)H-(v^x)^2-c^2 & (3-\gm)v^x & -(\gamma-1)v^y & \gm-1\\
      -v^xv^y & v^y & v^x & 0\\
      v^y[(\gm-2)H-c^2] & H-(\gm-1)(v^x)^2 & -(\gm-1)v^xv^y & \gm v^x
    \end{array}
  \right],
\end{equation}
where $c=\sqrt{\gm p/\rho}$ is the sound speed and $ H = E+\dfr p \rho = \dfr {(v^x)^2+(v^y)^2}2 +\dfr {c^2} {\gm-1}$ is the enthalpy. The Jacobian matrix $\pt\bbf/\pt\bu$ has four eigenvalues $ \lambda_1=v^x-c $, $ \lambda_2=\lambda_3=v^x $ and $ \lambda_4=v^x+c $ and four associated left eigenvectors
\begin{equation}\label{eq:lefteigen-Euler-f}
  \bga{l}
\bL_1^x = \dfr{1}{2c^2}\left[\dfr{\gamma-1}{2}((v^x)^2+(v^y)^2) + v^xc,\ -(\gamma-1)v^x-c,\ -(\gm-1)v^y,\ \gamma-1\right]^\top,\\[3mm]
\bL_2^x =-\dfr{1}{c^2}\left[\dfr{\gamma-1}{2}((v^x)^2+(v^y)^2) - c^2,\ -(\gamma-1)v^x,\ -(\gm-1)v^y,\ \gamma-1\right]^\top,\\[3mm]
\bL_3^x = \ \ \ \ \ \ \Big[-v^y,\ 0,\ 1,\ 0\Big]^\top,\\[3mm]
\bL_4^x = \dfr{1}{2c^2}\left[\dfr{\gamma-1}{2}((v^x)^2+(v^y)^2) - v^xc,\ -(\gamma-1)v^x+c,\ -(\gm-1)v^y,\ \gamma-1\right]^\top.
  \eda
\end{equation}

Along the boundary $ \Gm $, we have the solid-wall boundary condition $ v^x(0,y,t) = 0 $, which gives us $ u^{(0)}_2 = 0 $. Combining this with $ \bL_i^x\cdot\bu^{(0)} = w^{(0)}_i $, $ i=1,2,3 $, will specify the system \eqref{eq:general-0} of $ \bu^{(0)} $ for the current case as
\begin{equation}\label{eq:wall-2d-0}
\left\{
\bga{l}
  \bL_1^y\cdot\bu^{(0)} = w^{(0)}_1,\\[2mm]
  \bL_2^y\cdot\bu^{(0)} = w^{(0)}_2,\\[2mm]
  \bL_3^y\cdot\bu^{(0)} = w^{(0)}_3,\\[2mm]
  u^{(0)}_2 = 0.
\eda
\right.
\end{equation}
By taking temporal derivatives on both sides of $ u^{(0)}_2 = 0 $ and converting temporal derivatives to spatial ones,  one obtains
\begin{equation}\label{eq:wall-dx-dy}
[0, \ 1, \ 0, \ 0]\ \Big(-\dfr{\pt\bbf}{\pt\bu}(\bu^{(0)}) \ \dfr{\pt \bu}{\pt x}(0,y,t) -\dfr{\pt\bg}{\pt\bu}(\bu^{(0)}) \ \dfr{\pt \bu}{\pt y}(0,y,t) \Big)= 0,
\end{equation}
and therefore
\begin{equation}\label{eq:wall-dx-2d}
\bga{rl}
[0, \ 1, \ 0, \ 0]\ \dfr{\pt\bbf}{\pt\bu}(\bu^{(0)}) \ \bu^{(1)}= -[0, \ 0, \ 1, \ 0]\ \dfr{\pt\bg}{\pt\bu}(\bu^{(0)}) \ \dfr{\pt \bu}{\pt y}(0,y,t),
\eda
\end{equation}
where the Jacobian matrix of $\bg$ is 
\begin{equation}\label{eq:jacobi-Euler-g}
  \dfr{\pt \bg}{\pt \bu} = \left[
    \begin{array}{cccc}
      0&0&1&0\\
      -v^xv^y & v^y & v^x & 0\\
      (\gm-1)H-(v^y)^2-c^2 & -(\gamma-1)v^x & (3-\gm)v^y & \gm-1\\
      v^y[(\gm-2)H-c^2] & -(\gm-1)v^xv^y & H-(\gm-1)(v^y)^2 & \gm v^y
    \end{array}
  \right].
\end{equation}
Then the right hand side of \eqref{eq:wall-dx-2d} is calculated as 
\begin{equation*}
\bga{l}
\ \ \  -[0, \ 1, \ 0, \ 0] \ \dfr{\pt\bg}{\pt\bu}(\bu^{(0)}) \ \dfr{\pt \bu}{\pt y}(0,y,t)\\[2mm]
= [(v^x)^{(0)}(v^y)^{(0)},\ -(v^y)^{(0)},\ -(v^x)^{(0)},\ 0] \ \dfr{\pt \bu}{\pt y}(x,0,t)\\[2mm]
= (v^x)^{(0)}(v^y)^{(0)}\dfr{\pt u_1}{\pt y}(0,y,t) - (v^y)^{(0)}\dfr{\pt u_2}{\pt y}(0,y,t) - (v^x)^{(0)}\dfr{\pt u_3}{\pt y}(0,y,t)\\[2mm]
= 0. 
\eda
\end{equation*}
The last identity comes from the fact that $(v^x)^{(0)}=0$ and $\dfr{\pt u_2}{\pt y}(0,y,t)=0$. Hence we  combine \eqref{eq:wall-dx-2d} with $ \bL_i^y\cdot\bu^{(1)} = w^{(1)}_i $, $ i=1,2,3 $ to obtain
\begin{equation}\label{eq:wall-2d-1}
\left\{
\bga{l}
  \bL_1^y\cdot\bu^{(1)} = w^{(1)}_1,\\[2mm]
  \bL_2^y\cdot\bu^{(1)} = w^{(1)}_2,\\[2mm]
  \bL_3^y\cdot\bu^{(1)} = w^{(1)}_3,\\[2mm]
  \left[{ \dfr{\gm-1}{2}((v^y)^{(0)})^2,\ 0,\ -(\gamma-1)(v^y)^{(0)},\ \gm-1} \right] \bu^{(1)} = 0,
\eda
\right.
\end{equation}
which is the specific form of \eqref{eq:general-1} in the present case.

The two-dimensional numerical examples with the solid-wall boundary condition are presented in  Examples 7 and  8 in the next section when dealing with the double Mach reflection problem and the forward facing step problem.\\

\vspace{2mm}

\section{Numerical examples}\label{sec:numer}
 All examples displayed in this section are computed using the two-stage fourth-order scheme on Cartesian grids with the GRP solver \cite{Li-1} as the representative of the Lax--Wendroff type flow solvers. The fifth order Hermite-type WENO reconstruction is used for the spatial reconstruction. We denote this scheme by GRP4-HWENO5.\\

\noindent \textbf{Example 1. Linear scalar equations with smooth solutions. } We use a linear equation as the first example to verify the accuracy order of the current boundary condition treatment.  Consider the scalar IBVP \eqref{eq:IBVP} with the flux $ f(u)=u $. The initial and boundary conditions are $ u_0(x)=\sin(2\pi x) $ and $ g(t) = \sin(-2 \pi t) $, respectively. The setting of the initial data and the boundary condition allows the solution to be periodic. The inflow and outflow boundary condition treatments  are applied at $ x=0 $ and $ x=1 $, respectively.

The CFL number is set to be $ 0.4 $. The computation stops at $ t=5 $. The numerical errors and orders are shown in Table \ref{tab:scalar-advection}, which confirms  that the computation attains the expected order of accuracy.

\begin{table*}[!htbp]
  \centering
      \caption{The numerical errors and orders of the linear scalar equation in Example 1.}
  \begin{tabular}{|r|l|l|l|l|}
    \hline
$m$ & $ L_1 $ error&order&$ L_\infty $ error&order\\\hline
40  &    4.97e-07  &  4.32  &  1.75e-06  &  4.95\\
80  &    2.85e-08  &  4.13  &  8.64e-08  &  4.34\\
160  &   1.76e-09  &  4.02  &  5.66e-09  &  3.93\\
320  &   1.10e-10  &  4.00  &  3.60e-10  &  3.98\\
640  &   6.86e-12  &  4.00  &  2.27e-11  &  3.99\\\hline
  \end{tabular}
  \label{tab:scalar-advection}
\end{table*}
\vspace{2mm}

\noindent \textbf{Example 2. Nonlinear scalar equations.} This example purposes to see the accuracy order for  nonlinear equations. Consider the scalar IBVP \eqref{eq:IBVP} with the flux  $ f(u)={u^2}/{2} $ and the initial value $ u_0(x)=0.5+0.25\sin(2 \pi x) $. The boundary condition $ g(t) $ is given to be consistent with the initial value problem in which the initial data is periodically extend to $ x\in\Re $. At $ x=0 $, there are no explicit expressions available  for $ g(t) $ and its derivatives because of the nonlinearity. However,  the point-wise values of $ g $ could be obtained through the characteristic method and $g^{\prime\prime\prime}(t)$ used in \eqref{eq:inter-boundary-value} can be approximated by
\begin{equation}\label{eq:boundary-burgers}
\bga{l}
g^{\prime\prime\prime}(t) = \dfr{1}{2\tau^3}\big[-5g(t) + 18g({t+\tau}) - 24g(t+2\tau)\\
\ \ \ \ \ \ \ \ \ \ \ \ \ \ \ \ \ \ \ \ \ \ \ \ \ \ \ \ \ \ \ \ \ \ \ \ \ + 14g(t+3\tau) - 3g(t+4\tau)\big],
\eda
\end{equation}
where $\tau$ is proportional to the time step $k$. For example, we set $\tau=\dfr{k}{10}$ in our computations.

The CFL number is set to be $ 0.4 $. The computation stops at $ t=1/3\pi $.  For this example, $ x=0 $ is always an inflow boundary while $ x=1 $ is always an outflow boundary. Therefore the inflow and outflow boundary condition treatments  are applied at $ x=0 $ and $ x=1 $, respectively. The numerical errors and accuracy orders  in Table \ref{tab:scalar-burgers}  shows that the computation reaches the expected  order of accuracy.\\


\begin{table*}[!htbp]
  \centering
  \caption{The numerical errors and orders of the Burgers equation in Example 2.}
  \begin{tabular}{|r|l|l|l|l|}
    \hline
$m$ & $ L_1 $ error&order&$ L_\infty $ error&order\\\hline
40  &    6.00e-06  &  5.15  &  4.27e-05  &  4.07\\
80  &    1.49e-07  &  5.33  &  9.84e-07  &  5.44\\
160  &   7.63e-09  &  4.29  &  4.25e-08  &  4.53\\
320  &   4.77e-10  &  4.00  &  2.51e-09  &  4.09\\
640  &   2.95e-11  &  4.01  &  1.61e-10  &  3.96\\\hline

  \end{tabular}
  \label{tab:scalar-burgers}
\end{table*}
\vspace{2mm}

\noindent\textbf{Example 3. Scalar advection equations with discontinuous solutions.} In this example, we verify the capability of our boundary condition treatment to deal with the discontinuous boundary condition. The same example with similar initial and boundary conditions was used in \cite{ILW} for the same purpose.  Consider the IBVP \eqref{eq:IBVP} with the flux to be linear,  $ f(u) = u $. The initial data is $ u_0(x)=\sin(4\pi x) $ and the boundary data is taken as
\begin{equation}
g(t) = \left\{
\begin{array}{lr}
  \sin(-4\pi t),& t < 0.25,\\
  0,& 0.25 < t < 0.5,\\
  3, &t > 0.5.\\
\end{array}
\right.
\end{equation}

\begin{figure}[!htb]
\centering
 \includegraphics[width=.8\textwidth]{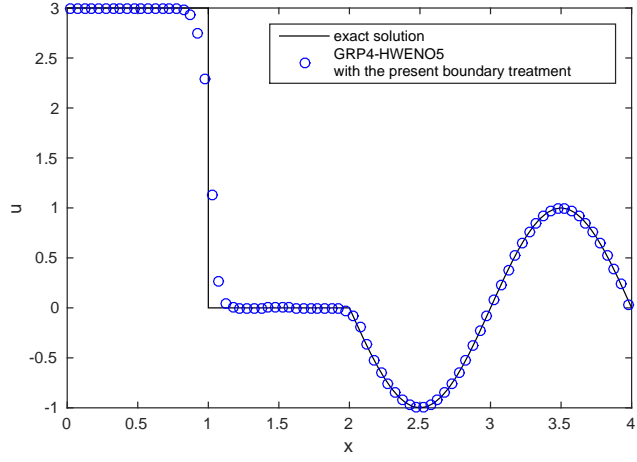}
 \caption{The discontinuous solution in Example 3 with $80$ computational cells (circles). The exact solution  is shown as a solid curve.}
 \label{fig:scalar-disc}
\end{figure}

The inflow and outflow boundary condition treatments are applied at $ x=0 $ and $ x=1 $, respectively. We compare the result in the domain $ (0,1) $ with the exact one in Figure \ref{fig:scalar-disc}. The numerical solution matches the exact solution represented by the solid line perfectly.\\

\vspace{2mm}

\noindent \textbf{Example 4. Linear systems,} Consider the linearized Euler equations
\begin{equation}
\label{linear-Euler}
\dfr{\pt}{\pt t}\left[
\begin{array}{c}
\hat{\rho}\\
\hat{v}\\
\hat{p}
\end{array}
\right]+
\left[
\begin{array}{ccc}
    \bar{v}&\bar{\rho}&0\\
    0&\bar{v}&{1}/{\bar{\rho}}\\
    0&\gm\bar{p}&\bar{v}
\end{array}
\right] \dfr{\pt}{\pt x}\left[
\begin{array}{c}
\hat{\rho}\\
\hat{v}\\
\hat{p}
\end{array}
\right]=0,\ \ \  x\in(0,2),\ \  \  t>0
\end{equation}
with the background state $ (\bar{\rho},\bar{v},\bar{p}) $ and the perturbation $ (\hat{\rho},\hat{v},\hat{p}) $, where $ \gm=1.4 $ is the specific heat ratio. The eigenvalues of $ \bar{A} $ are $ \lambda_1= \bar{v}-\bar{c}$, $ \lambda_2= \bar{v}$ and $ \lambda_3= \bar{v}+\bar{c}$ with $ \bar{c}^2=\gm\bar{p}/\bar{\rho}$. In this example, we set $ \bar{\rho}=1.4 $, $ \bar{v}=1 $ and $ \bar{p}=4 $, which makes $ \lambda_1 < 0 < \lambda_2 < \lambda_3 $. Therefore subsonic inflow and subsonic outflow boundary conditions are prescribed at $ x=0 $ and $ x=2 $, respectively.

This example is designed to allow the solution to be a combination of three sine waves carried by the three characteristics, i.e.
\begin{equation}\label{eq:sol-linear-sys}
\left[
\begin{array}{c}
\hat{\rho}(x,t)\\
\hat{v}(x,t)\\
\hat{p}(x,t)
\end{array}
\right] = \Big(\bR_1, \bR_2, \bR_3\Big)\left[
\begin{array}{l}
\alpha_1\sin(k_1\pi(x-(\bar{v}-\bar{c})t))\\
\alpha_2\sin(k_2\pi(x-\bar{v}t))\\
\alpha_3\sin(k_3\pi(x-(\bar{v}+\bar{c})t))
\end{array}
\right],
\end{equation}
where $ k_i $ and $ \alpha_i $, $ i=1,2,3 $ are parameters, and
\begin{equation}
  \bR_1=\Big[
      1, \ -\dfr{\bar{c}}{\bar{\rho}}, \ \bar{c}^2
  \Big]^\top, \ \ \bR_2=\Big[
      1, \ 0, \ 0
  \Big]^\top, \ \ \bR_3=\Big[
      1, \ \dfr{\bar{c}}{\bar{\rho}}, \ \bar{c}^2
  \Big]^\top.
\end{equation}
For example, we set $ k_1=1 $, $ k_2=3 $, $ k_3=2 $, $ \alpha_1=0.1 $, $ \alpha_2=-0.1 $ and $ \alpha_3=0.08 $. Following the instruction of the fluid dynamics \cite{gas-dynamics}, the inflow density $\rho_\text{in}(t)=\hat{\rho}(0,t)$ and pressure $p_\text{in}(t)=\hat{p}(0,t)$ are given at $ x=0 $ and the outflow pressure $p_\text{ex}(t)=\hat{p}(2,t)$ is given at $ x=2 $. The initial condition is defined by setting $ t=0 $.

In the computation, subsonic inflow and subsonic outflow boundary condition treatments are applied at $ x=0 $ and $ x=2 $, respectively. The CFL condition is $ 0.4 $ and the output time is $ t=10 $. Table \ref{tab:system-advection} shows the numerical errors and orders.

\begin{table*}[!htbp]
  \centering
  \caption{The numerical errors and orders of the linear system in Example 4.}
  \begin{tabular}{|r|l|l|l|l|}
    \hline
$m$ & $ L_1 $ error&order&$ L_\infty $ error&order\\\hline
40  &  1.87e-2  &  -0.12   &  2.11e-2  &  0.06 \\
80  &  8.18e-4  &  4.52   &  1.01e-3  &  4.39 \\
160  &  4.08e-5  &  4.33   &  4.54e-5  &  4.47 \\
320  &  2.79e-6  &  3.87   &  3.04e-6  &  3.90 \\
640  &  1.78e-7  &  3.97   &  1.89e-7  &  4.00 \\\hline
  \end{tabular}
  \label{tab:system-advection}
\end{table*}
\vspace{2mm}

\noindent \textbf{Example 5. The Woodward--Colella problem.}  This classical example  assumes that initially the gas is at rest and ideal with $\gm=1.4$ in the computation domain $ [0,1] $, the density is everywhere unit and the pressure is $p=1000$ for $0\leq x <0.1$ and $p=100$ for $0.9<x\leq 1.0$, while it is only $p=0.01$ for $0.1<x<0.9$. The solid-wall boundary condition is prescribed at both ends. Two numerical methods are used to deal with the  boundary condition. The first one is the traditional numerical treatment with which we  symmetrically extend the solution values into ghost cells. The second one uses the present boundary condition treatment in the solid-wall case, with which the procedure in  Subsections \ref{subsec:system-general} and \ref{subsec:wall} are applied. The CFL number is $ 0.6 $.

\begin{figure}[!htb]
\centering
 \includegraphics[width=\textwidth]{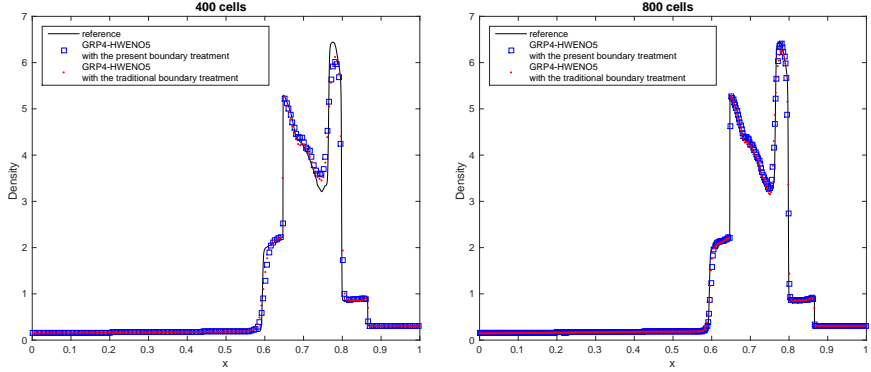}
 \caption{The Woodward--Colella problem computed with the present boundary condition treatment (squares) and the conventional reflection boundary condition treatment (dots) with 400 cells (200 are shown, left) and 800 cells (400 are shown, right). The numerical scheme used is the GRP4-HWENO5 scheme. The solid lines are the reference solution computed with 4000 cells.}
 \label{fig:wood}
\end{figure}

The numerical solutions using the present boundary condition treatments are displayed in Figure \ref{fig:wood} at output time $t=0.038$, in comparison with those using symmetric extension of interior information to ghost cells. The similar results verifies the effectiveness of the current approach for the one-dimensional solid-wall boundary condition.

\vspace{2mm}

\noindent \textbf{Example 6. The nozzle flow.}  The problem of the nozzle flow is  quasi one-dimensional. A converging-diverging duct occupies the spatial interval $ x \in  (0,1) $ and has a continuous cross-sectional area function $ A(x) $ given by
\begin{equation}\label{eq:nozzle-A}
  A(x) =
\left\{
\bga{lr}
A_{\text{in}}\exp[-\log(A_{\text{in}})\sin^2(2\pi x)], & 0 \leq x\leq 0.25,\\[2mm]
A_{\text{ex}}\exp[-\log(A_{\text{ex}})\sin^2( \dfr{2\pi(1-x)}{3})], & 0.25 \leq x\leq 1,
\eda
\right.
\end{equation}
with $A_{\text{in}}=4.864317646$ and $A_{\text{ex}}=4.234567901$. The cross-sectional area reaches its minimal value at $ x=0.25 $, which is called the throat of the duct. The governing PDEs of the nozzle flow are the Euler equations with the geometric source term \eqref{eq:nozzle}. The fluid flows from the left to the right.  We set $ x = 0 $ as  the entrance of the duct and $ x = 1 $ as the exit. The flow in the duct should finally reaches a steady state as the physics indicates.

There are two types of steady states: a continuous steady state and a discontinuous steady state containing a stationary shock. The initial conditions for both cases take
\begin{equation}
\label{eq:init-nozzle}
  (\rho(x,0),\ v(x,0),\ p(x,0))=\left\{
\bga{lr}
(\rho_0, \ 0, \ p_0), & x < 0.25,\\[2mm]
(\rho_0 (p_{\text{ex}} / p_0)^{1/\gamma}, \ 0, \ p_{\text{ex}}), & \ \ \ x>0.25,
\eda
\right.
\end{equation}
where $\gm=1.4$, $\rho_0$ and $p_0$ are parameters, determining if the steady state is continuous or not. In the previous numerical studies of the nozzle flow \cite{Ben-Artzi-01,Li-2},   the inflow density, velocity and pressure are assigned as the inflow boundary condition   to the ghost cells out of the entrance, and  the outflow pressure is assigned as the outflow boundary condition to the ghost cell out of the exit. In the present study, the approximation  strategy of boundary conditions  in Subsections \ref{subsec:system-general} and \ref{subsec:in-and-out} is  applied.

\begin{figure}[!htb]
\centering
 \includegraphics[width=\textwidth]{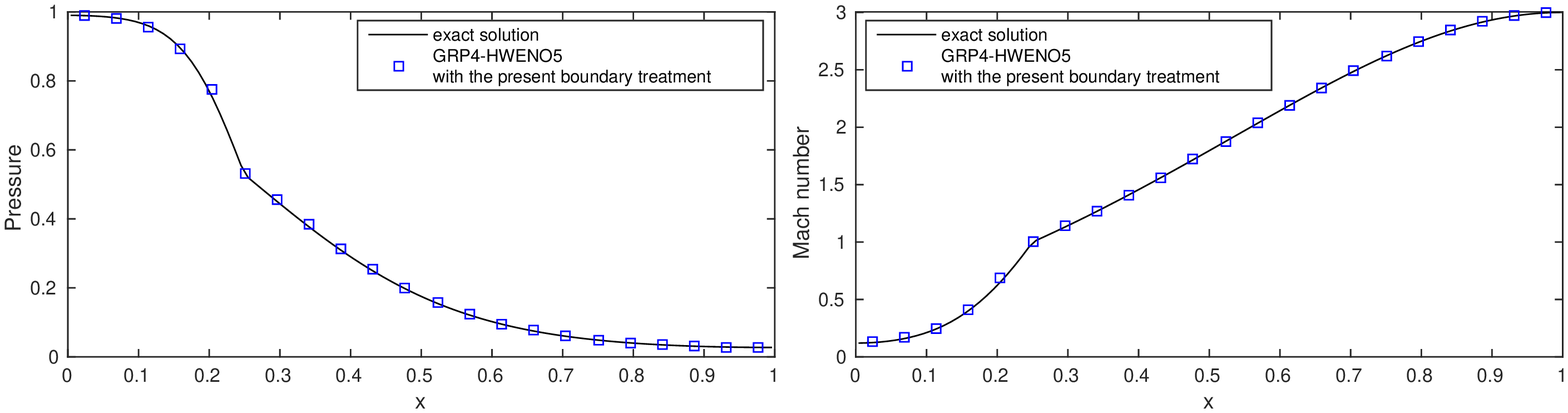}
 \caption{The isentropic flow throughout all the duct computed with the two-stage fourth-order scheme. The density and the Mach number at $t=5$ are shown (squares). 22 cells are used. The solid line is the exact solution given by \eqref{eq:isen-nozzle}.}
 \label{fig:nozzle-0}
\end{figure}

For the first case, we set $ \rho_0 = p_0 = 1 $ and $ p_\text{ex}  = 0.0272237$ in \eqref{eq:init-nozzle}. See \cite{Ben-Artzi-01}. This makes the steady solution a continuous isentropic one defined by
\begin{equation}\label{eq:isen-nozzle}
\bga{l}
  \rho(x) = \rho_0\left(1+\dfr{\gamma-1}{2}[M(x)]^2\right)^{-\frac{1}{\gamma-1}}, \\[3mm]
  p(x) = p_0\left(1+\dfr{\gamma-1}{2}[M(x)]^2\right)^{-\frac{\gamma}{\gamma-1}},\\[4mm]
  v(x) = M(x)\sqrt{\gamma \ p(x)/\rho(x)},
\eda
\end{equation}
where $M(x)$ is determined by the sectional area $A(x)$ through the relation
\begin{equation}\label{eq:mach}
  [A(x)]^2 = \dfr{1}{[M(x)]^2}\left[\dfr{1}{\gamma+1}\left(1+\dfr{\gamma-1}{2}[M(x)]^2\right)\right]^\frac{\gamma+1}{\gamma-1}.
\end{equation}
In this case, the flow is subsonic upstream to the throat and supersonic downstream to the throat. Thus  the inflow boundary condition at the entrance $x=0$ is prescribed as in  \eqref{eq:nozzle-in} with,
\begin{equation}\label{eq:nozzle-in-para}
\bga{l}
  p_\text{in} := p_0\left(1+\dfr{\gamma-1}{2}[M(0)]^2\right)^{-\frac{\gamma}{\gamma-1}},\\
  \rho_\text{in} := \rho_0\left(1+\dfr{\gamma-1}{2}[M(0)]^2\right)^{-\frac{1}{\gamma-1}}.
\eda
\end{equation}
At the exit  $x=1$, the flow  is supersonic and no boundary condition is required. The numerical result is displayed in Figure \ref{fig:nozzle-0} with the current method, using 22 computational cells. The output time is $t=5$ and the CFL number is $0.6$. The solution converges to the expected steady one and attains a better agreement with the steady solution compared with those given in \cite{Ben-Artzi-01,Li-2}.

\begin{figure}[!htb]
\centering
 \includegraphics[width=\textwidth]{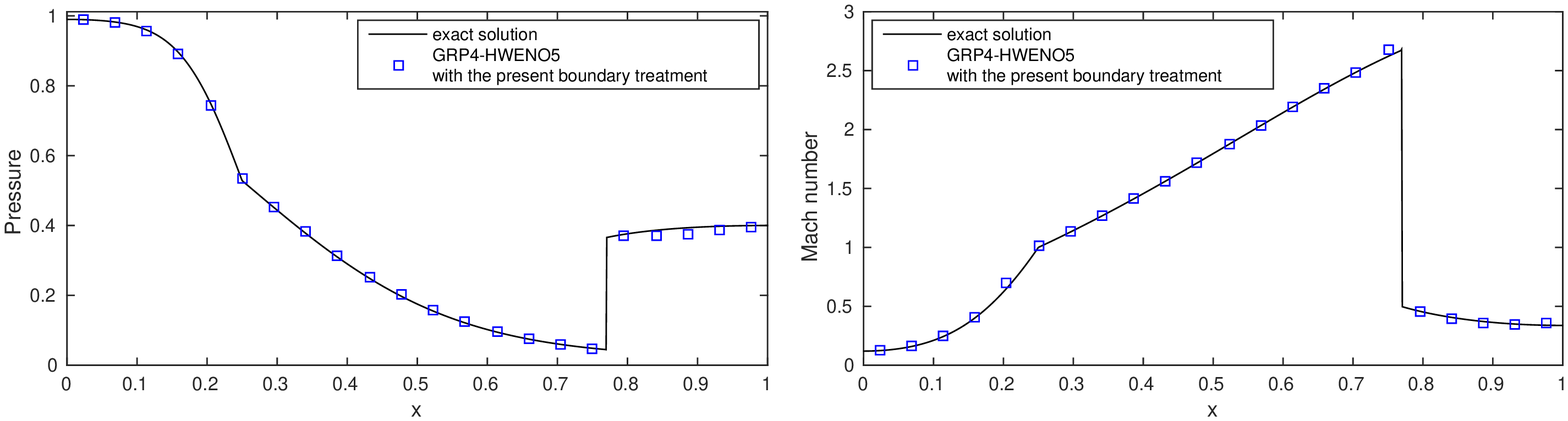}
 \caption{The flow with a steady shock computed with the two-stage fourth-order scheme. The density and the Mach number at $t=5$ are shown (squares). 22 cells are used. The solid line is the exact solution given by \eqref{eq:isen-nozzle}.}
 \label{fig:nozzle-1}
\end{figure}

The other steady solution  contains a stationary shock separating two pieces of isentropic solutions defined as in \eqref{eq:isen-nozzle} with separate pairs of $(\rho_0,p_0)$. In this case we set $\rho_0=p_0=1$ and $p_\text{ex}=0.4$ in \eqref{eq:init-nozzle} to get the initial data.
The shock stands downstream to the throat, and the flow jumps from supersonic to subsonic after passing the shock. The exit $ x=1 $ is subsonic for such a  case. Now we set the inflow boundary condition to be \eqref{eq:nozzle-in} and \eqref{eq:nozzle-in-para} with $\rho_0=p_0=1$ at the entrance $x=0$ and the outflow boundary condition to be \eqref{eq:nozzle-ex} with $p_\text{ex}=0.4$ at the exit $x=1$.  Figure \ref{fig:nozzle-1} shows the numerical results with 22 computational cells. The output time is $t=5$ and the CFL number is $ 0.6 $. Once again, the solution converges to the expected steady one  and matches it better compared with those given in \cite{Ben-Artzi-01,Li-2}.

As the accuracy test is performed for such a case, we need to modify the cross section a little bit, such as
\begin{equation}\label{eq:nozzle-A2}
  A(x) = A_{\text{in}}\exp[-\log(A_{\text{in}})\sin^2(\pi x)], \ \ \text{for} \ \ 0 \leq x\leq 1,
\end{equation}
instead of \eqref{eq:nozzle-A}, in order to guarantee the flow is smooth. This is because $A''(x)$ in \eqref{eq:nozzle-A} is discontinuous, which leads to the solution formula in  \eqref{eq:isen-nozzle} has the discontinuity in its  first order derivative at the throat of the duct which can be seen in Figure \ref{fig:nozzle-0}. 
The initial data in this case is
\begin{equation*}
\label{eq:init-nozzle-A2}
  (\rho(x,0),\ v(x,0),\ p(x,0))=\left\{
\bga{lr}
(\rho_0, \ 0, \ p_0), & x < 0.5,\\[2mm]
(\rho_0 (p_{\text{ex}} / p_0)^{1/\gamma}, \ 0, \ p_{\text{ex}}), & \ \ \ x>0.5,
\eda
\right.
\end{equation*}
by setting $\rho_0=p_0=1$ and $p_\text{ex}=0.021910717$. The numerical errors and accuracy orders of the momentum $A\rho v$ are shown in Table \ref{tab:accuracy-nozzle-smooth},  which verifies the numerical accuracy of the present boundary condition treatment.
\begin{table*}[!htbp]
  \centering
  \caption{Numerical errors and accuracy orders for the case  with the cross section \eqref{eq:nozzle-A2}}
  \begin{tabular}{|r|l|l|l|l|}
    \hline
$m$ & $ L_1 $ error&order&$ L_\infty $ error&order\\\hline
 40  &  2.13e-04  &  2.04  &  2.59e-03  &  1.94\\
 80  &  1.40e-06  &  7.25  &  1.22e-05  &  7.73\\
160  &  2.87e-08  &  5.61  &  3.34e-07  &  5.19\\
320  &  2.10e-09  &  3.77  &  2.12e-08  &  3.98\\
640  &  1.30e-10  &  4.02  &  1.34e-09  &  3.98\\ \hline
  \end{tabular}
  \label{tab:accuracy-nozzle-smooth}
\end{table*}

\vspace{2mm}

\noindent \textbf{Example 7. The double Mach reflection problem.}  This is a standard two-dimensional test problem for high resolution schemes. The computational domain for this problem is $ [0,4] \times [0,1] $, and $ [0,3] \times [0,1]$ is shown. A reflective wall lies at the bottom of the computational domain starting from $ x=\frac{1}{6} $. Initially a right-moving Mach 10 shock is positioned at $ x=\frac{1}{6} ,y=0$ and makes a $ \frac{\pi}{3} $ angle with the $x$-axis. More details about this problem can be seen in \cite{Woodward-Colella}.

Our computations use both the traditional and the present new numerical boundary condition treatments to deal with the reflective boundary condition along the wall $ \{(x,0): x\in[1/6, 4]\} $. For the traditional boundary condition treatment, we use the symmetrical extension for  the ghost cells outside the boundary. For the present new boundary condition treatment, we apply the procedure in Subsections \ref{subsec:system-general} and \ref{subsec:wall-2d}. 

\begin{figure}[!htb]
\centering
 \includegraphics[width=\textwidth]{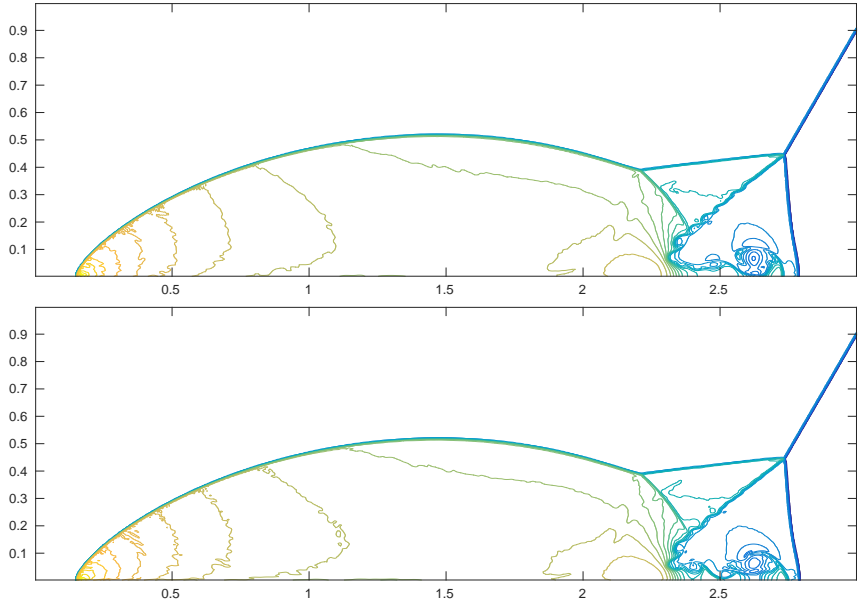}
 \caption{The numerical results of the double mach reflection problem in Example 7 given by the scheme GRP4-HWENO5 combined with the traditional reflection boundary condition treatment (upper) and the current boundary condition treatment (lower). The contours of density are shown. $ 960\times 240 $ cells are used.}
 \label{fig:DM-240}
\end{figure}

The computations are carried out by the GRP4-HWENO5 scheme with $960\times 240$ grids, using both numerical boundary condition treatments. The numerical results are displayed in Figure  \ref{fig:DM-240}  with $ 30 $ contours of the density at time $ 0.2 $. The CFL number is $ 0.6 $. The two numerical results are similar which indicates that the present approach works well for the two-dimensional solid-wall boundary condition.

 \vspace{0.2cm}

\noindent\textbf{Example 8. The forward facing step problem.} This is also a standard test problem for the two-dimensional computations proposed in \cite{Woodward-Colella}. The wind tunnel is $1$ length unit wide and $3$ length units long. The step is $0.2$ length units high and is located $0.6$ length units from the left-hand end of the tunnel. The problem is initialized by a unit right-going Mach $3$ flow with $ (\rho_0, u_0, v_0, p_0)=(1.4,3,0,1) $ in the tunnel. Reflective boundary conditions are applied along the walls of the tunnel.

\begin{figure}[!htb]
\centering
 \includegraphics[width=\textwidth]{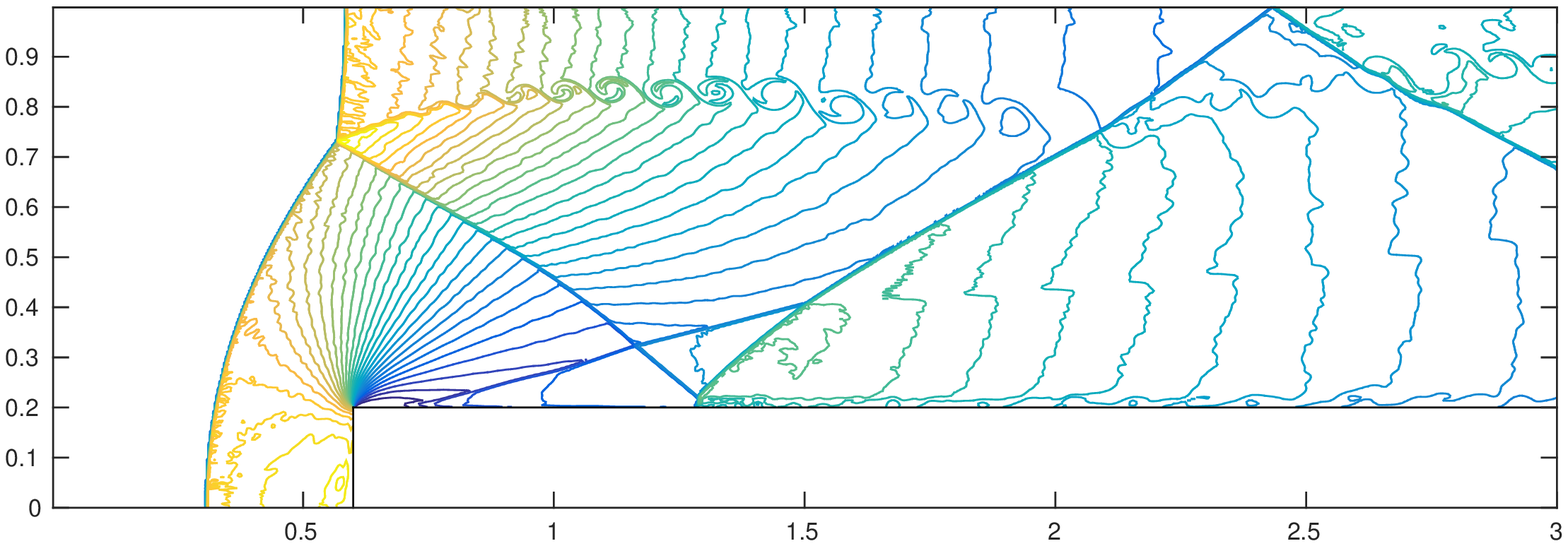}
 \includegraphics[width=\textwidth]{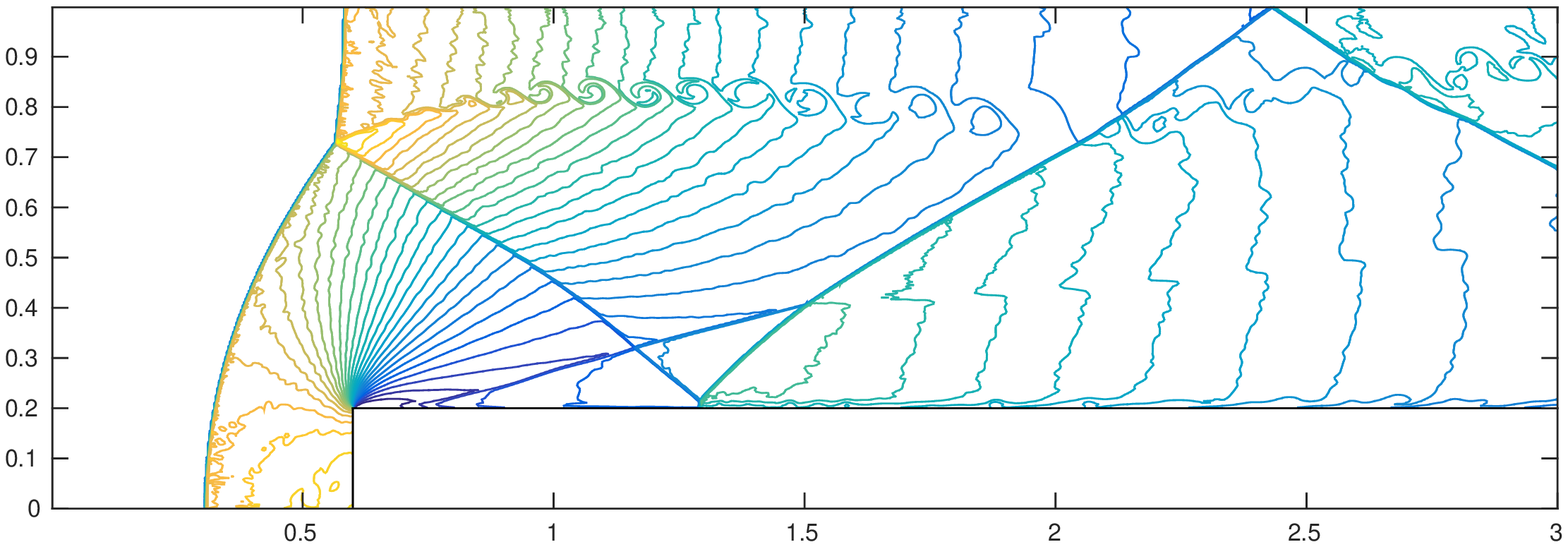}
 \caption{ The numerical results of the forward facing step problem in Example 8 given by the scheme GRP4-HWENO5 using the traditional boundary condition treatment (upper) and the present one (lower). The contours of density are shown. $ 960\times 320 $ cells are used.}
 \label{fig:FS-NBT-2}
\end{figure}

As in Example 7, both the traditional and the present new boundary condition treatments are applied to the reflective walls of the tunnel. 
The numerical results are shown  in Figure \ref{fig:FS-NBT-2}, with  $ 960\times 320 $ cells. The computations stop at the time $ t=4 $ and the CFL number is $ 0.6 $. The two numerical results are similar which once again verifies the effectiveness of the present approach.

\vspace{2mm}

\section{Discussions}

In this paper we provide a method to approximate boundary conditions with fourth-order accuracy   in order to suit the two-stage fourth order accurate schemes for hyperbolic conservation laws that we proposed earlier \cite{Du-Li-1}. The application is specified to the compressible Euler equations with several commonly-used boundary conditions. We are certainly aware that there are many issues waiting for investigation, such as moving boundary problems, solid wall boundary conditions with curved geometry, small cut-cell problems, and  moving boundary problems,  which will be studied in a future paper.  \\

Here we would like to emphasize that the inverse Lax-Wendroff approach is fundamental in the sense that the governing equations are effectively adopted to treat the boundary conditions \cite{ILW-HJ-1,ILW,high-order-ILW,ILW-stb}. In the paper the inverse Lax-Wendroff approach is used with  the least complexity  by taking advantage of the two-stage method \cite{Du-Li-1}. No successive differentiation of governing equations is made.
Nevertheless, the boundary conditions have to be carefully treated at intermediate stages, analogous to any other multi-stage numerical methods.  \\

Of course, just as the fact that the two-stage fourth order method in \cite{Du-Li-1} can be extended to higher order accurate multi-stage method in \cite{Xu-Li}, the current numerical boundary condition treatment can be extended to much higher order accuracy in a straightforward way.  The application to hyperbolic systems beyond the compressible fluid flows  can be treated similarly.  We do not want to repeat the technicality from the scientific viewpoint.

\vspace{2mm}

\appendix
\section{The interpolation results in subsection \ref{sec:inflow}}
This appendix is dedicated to list the interpolation results  in Section \ref{sec:inflow}. Recall that we assume $ x=0 $ and $ x=1 $ are the inflow and outflow boundaries for the IBVP \eqref{eq:IBVP} of the one-dimensional scalar conservation laws, respectively. The stencils are denoted in \eqref{stencil}. \\

\subsection{Cell averages and cell differences}\label{app:low-order}
The reconstructed average of $ u $ in $ I_{-1} $ and $ I_{-2} $ on those stencils are:
\begin{equation}
\bga{l}
      \bar{u}_{-1}^{(2)} = \dfr 14 (-6g + 6~h~ f^\prime(g)^{-1} \ g^{\prime} + 11\bar{u}_{0} - \bar{u}_{1}),\\[2mm]
      \bar{u}_{-1}^{(1)} = h~ f^\prime(g)^{-1} \ g^{\prime} + \bar{u}_{0},\\[2mm]
      \bar{u}_{-1}^{(0)} = g + \dfr 12~ h~ f^\prime(g)^{-1} \ g^{\prime},\\[4mm]
      \bar{u}_{-2}^{(2)} = \dfr 14 (-90g + 42~h~ f^\prime(g)^{-1} \ g^{\prime} + 105\bar{u}_{0} - 11\bar{u}_{1}),\\[2mm]
      \bar{u}_{-2}^{(1)} = -6gh + 5~f^\prime(g)^{-1} \ g^{\prime} + 7\bar{u}_{0},\\[2mm]
      \bar{u}_{-2}^{(0)} = g + \dfr 32~ h~ f^\prime(g)^{-1} \ g^{\prime}.
\eda
\end{equation}
The reconstructed $ x $-difference of $ u $ in $ I_{-1} $ and $ I_{-2} $ on those stencils are:
\begin{equation}
\bga{l}
      \Delta u_{-1}^{(2)} = \dfr 1{8h} (66g - 34~h~ f^\prime(g)^{-1} \ g^{\prime} - 73\bar{u}_{0} + 7\bar{u}_{1}),\\[2mm]
      \Delta u_{-1}^{(1)} = \dfr 1{2h} (6g - 5~h~ f^\prime(g)^{-1} \ g^{\prime} - 6\bar{u}_{0}),\\[2mm]
      \Delta u_{-1}^{(0)} = -f^\prime(g)^{-1} \ g^{\prime},\\[4mm]

      \Delta u_{-2}^{(2)} = \dfr 1{8h} (294g - 118~h~ f^\prime(g)^{-1} \ g^{\prime} - 331\bar{u}_{0} + 37\bar{u}_{1}),\\[2mm]
      \Delta u_{-2}^{(1)} = \dfr 1{2h} (18g - 11~h~ f^\prime(g)^{-1} \ g^{\prime} - 18\bar{u}_{0}),\\[2mm]
      \Delta u_{-2}^{(0)} = -f^\prime(g)^{-1} \ g^{\prime}.
\eda
\end{equation}

\subsection{Smoothness indicators}\label{app:SI}
The smoothness indicators on these stencils are listed as follows,
\begin{equation}
\bga{l}
\beta^{(2)} = \dfr 1{80}\Big[    66516 g^2 +     9444 (h f^\prime(g)^{-1}g^\prime)^2 - 56348  f^\prime(g)^{-1}g^\prime h \bar{ u}_0\\[2mm]
\ \ \ \ \ \ \ \ \ \ \ \ \ \ \ \ \ + 85929 \bar{ u}_0^2 + 6644  f^\prime(g)^{-1}g^\prime h \bar{ u}_1 -     20694 \bar{ u}_0 \bar{ u}_1 +     1281 \bar{ u}_1^2\\[2mm]
\ \ \ \ \ \ \ \ \ \ \ \ \ \ \ \ \ + 12 g (    4142  f^\prime(g)^{-1}g^\prime h -    12597 \bar{ u}_0 +     1511 \bar{ u}_1)\Big],\\[3mm]

\beta^{(1)} = 48 g^2 + 54 g h  f^\prime(g)^{-1}g^\prime + 16 (h f^\prime(g)^{-1}g^\prime)^2\\[2mm]
\ \ \ \ \ \ \ \ \ \ \ \ \ - 96 g \bar{ u}_0 + 48 \bar{ u}_0^2 - 54 h f^\prime(g)^{-1}g^\prime \bar{ u}_0,\\[3mm]

\beta^{(0)} = (h f^\prime(g)^{-1}g^\prime)^2.
\eda
\end{equation}

\section*{Acknowledgements}
This work is supported by NSFC (no. 11771054) and Foundation of LCP.  Both authors appreciate an anonymous reviewer for his suggestions on the accuracy test of the nozzle problem and  on curved geometry or small cut-cell problems.  The latter suggestion would be adopted in the future due to very technical details and the clarity of the  paper. We also thank the other two reviewers for their  suggestions that  substantially polish this paper.

\vspace{0.5cm}
\end{document}